\documentclass {amsart}

\newtheorem{theorem}{Theorem}[section]
\newtheorem{lemma}[theorem]{Lemma}
\newtheorem{proposition}[theorem]{Proposition}
\newtheorem{corollary}[theorem]{Corollary}

\theoremstyle{definition}

\theoremstyle{remark}
\newtheorem{remark}[theorem]{Remark}

\numberwithin{equation}{section}

\usepackage{amssymb}
\usepackage{amscd}
\usepackage{amsmath}
\usepackage{graphicx}

\usepackage[all]{xy}

\newcommand{\ff}{\mathbb F}

\newcommand{\pp}{\mathbb P}
\newcommand{\calO}{\mathcal{O}}
\newcommand{\calE}{\mathcal{E}}

\newcommand{\st}{\overline}

\newenvironment{dimi}{\small}

\begin{document}

\title[On the classification of surfaces of General type]{On the classification of surfaces of general type
 with non-birational bicanonical map and Du Val Double planes}
\author{Giuseppe Borrelli}
\address{}
\curraddr{}
\email{borrelli@mat.uniroma3.it}
\thanks{This work was partially sopported by EU Research Training Network EAGER (HPRN-CT-2000-00099).}


\date{}
\dedicatory{This paper is dedicated to the memory of prof. Paolo Francia.}


\begin{abstract}
We classify surfaces of general type whose bicanonical map $\varphi \sb{2K}$ is composed with
a rational map of degree 2 onto a rational or ruled surface.
In particular, this is always the case if $q=0,p_g\ge2 $ and $\varphi_{2K}$ is not birational.

We prove that such a surface $S$ either has a genus 2 pencil or is the smooth model of a double plane branched along a
reduced curve with certain singularities, a configuration already suggested by Du Val in the 1950's.

In the last case we show that $S$ has a rational pencil $|C|$ such that the general member is a smooth hyperelliptic curve
of genus $3$, unless $K_S$ is ample and either $p_g(S)=6,K_S^2=8$ or $p_g(S)=3,K_S^2=2$.
\end{abstract}

\maketitle

Let $S$ be a smooth minimal algebraic surface over the complex numbers with geometric genus $p_g(S)=h^0(S,\calO_S(K_S))$
and irregularity $q(S)=h^1(S,\calO_S)$. Assume that $S$
is of general type, then the bicanonical map of $S$ is the rational map
\[
\varphi_{2K}:S\dashrightarrow S_2\subseteq \pp^{K_S^2+p_g(S)-q(S)}
\]
defined by the linear system $|2K_S|$, 
where $K_S$ is a canonical divisor for $S$ and $S_2$ is the bicanonical image.

A theorem of Xiao \cite{X1} says that $S_2$ is a surface unless $p_g(S)=q(S)=0$ and $K_S^2=1$. On the other hand
there is a $standard$ $case$ for the non birationality of $\varphi_{2K}$, that is if $S$ has a pencil
$|C|$ such that the general element $C\in |C|$ is a curve of genus 2.

By \cite{R} if $\varphi_{2K}$ is not birational and $S$ does not present the standard case then $K_S^2\leq 9$,
thus there are finitely many
families of such surfaces and it is natural to study and try to classify them.

In the 1950's Du Val suggested that examples of minimal surfaces of general type
with non birational bicanonical map can be obtained in the following way.

Let $X$ be a smooth surface and $G\subset X$ a reduced curve such that \textit{
\begin{itemize}
\item [$\mathcal{B}$)] either $X=\ff_2$ and $G=C_0+G^\prime$, where $G^\prime\in |7C_0+14\Gamma|$ and $G^\prime$
has at most non essential singularities;
\item [$\mathcal{D}$)] or $X=\pp^2$ and $G$ is a smooth curve of degree $8$;
\item [$\mathcal{D}_n$)] or $X=\pp^2$ and $G=G'+L_1+\dots +L_n$, with $ n\in \{0,1,\dots,6\}$ $(G=G^\prime$ if $n=0)$,
where $L_1,\dots,L_n$ are distinct lines meeting at a point $\gamma$ and $G'$ is a curve of \hbox {degree $10+n$}.
The singularities of $G$, besides the non essential ones, are a $(2n+2)$-tuple point at $\gamma$,
a $[5,5]$-point lying on $L_i$, $i=1,\dots,n$,
possibly some $4$-tuple points or $[3,3]$-points;
\end{itemize}
}
then $S$ is the smooth minimal model of the double cover $X^\prime \rightarrow X$ branched along $G$. Here $\ff_2$ is the
Hirzebruch surface $\textbf{P}(\calO_{\pp^1}\oplus\calO_{\pp^1}(2))$ and $\Gamma,C_0$ its fibre and negative section with $C_0^2=-2$.

We will refer to such examples as the $Du$ $Val$ $examples$, whilst by abuse of notation we will say that
$X^\prime$ is a $Du$ $Val$ $double $ $plane$ (of type $\mathcal{B},\mathcal {D}$ or $\mathcal{D}_n$ respectively)

Under the hypotheses $h^1(S,\calO_S)=0, p_g(S)\ge 3$ and that the general canonical curve is irreducible 
Du Val proved that if $\varphi_{2K}$ is not birational and $S$ does not present the standard case then
$S$ is one of the above examples.

More recently C.Ciliberto, P.Francia and M.Mendes Lopes have considered the same problem
in \cite{CFM} and \cite{CM} removing the hypothesis
concerning the general canonical curve and the regularity of $S$.
They worked it out with modern arguments and essentially they confirmed
the classification of Du Val for the regular case (i.e. $q(S)=0$).

In my PhD thesis (cfr. \cite{B}) I proved an analogous result
for regular surfaces with $p_g(S)=2$ under the assumption that the canonical system has no fixed part.

In this article we extend the above results rephrasing Du Val's claim.
For this we remark that if $q(S)=0$ and $p_g(S)\ge 2$ then
$\varphi_{2K}$ is either birational or a (generically finite) morphism of degree 2 onto a rational surface.

In fact, $\varphi_{2K}$ has no base points by \cite{F} and writing
$|K_S|=|M|+F$ where $|M|$ is the movable part we have that the general curve $M\in |M|$ is irreducible and
$|2K_S|$ separates different curves of $|M|$.
Therefore, looking at the exact sequence
\[
H^0(S,\calO_S(K_S+M))\rightarrow H^0(M,\calO_M(K_M))\rightarrow 0
\]
we get that if $\varphi_{2K}$ is not birational the rational map $\varphi_{|K_S+M|}$ defined by the linear system
$|K_S+M|\subset|2K_S|$ is not birational on a general $M$.
Hence $M$ is hyperelliptic and $\varphi_{2K}:S\rightarrow S_2$ is a generically finite morphism of \hbox{degree 2}.
Therefore, $S_2$ is a surface of degree $2K_S^2$ in $\pp^N$ where $N=K_S^2+p_g(S)$ and as 
$2K_S^2<2N-2$, $S_2$ is a ruled surface. Whence, $S_2$ is rational since $S$ is regular.

More generally, we may consider minimal surfaces of general type for which
the bicanonical map factors through a rational map $\phi$ of degree 2
onto a rational or ruled surface, that is if there exists a commutative diagram
\[
 \xymatrix@!0{
  S \ar@{.>} [dd]_{\phi} \ar@{.>}[rrr]^{\varphi_{2K}} & && S_2    \\
 \\
\Sigma \ar@{.>}[uurrr]_{\phi_2}
   }
\]
where $\phi $ is a (generically finite) rational map of degree two and $\Sigma$ is a rational or ruled surface.

Our main result is the following

\begin{theorem}\label{mainth}
Let $S$ be a smooth minimal surface of general type  which does not present the standard case.
Then the following three conditions are equivalent:
\begin{itemize}
\item[$a)$] the bicanonical map of $S$
factors through a rational map of degree 2 onto a rational or ruled surface
\item[$b$)] the bicanonical map of $S$
factors through a rational map of degree 2 onto a rational surface
\item[$c$)] $S$ is the smooth minimal model of a Du Val double plane.
\end{itemize}
Moreover, let $S$ be as in (c) (resp. (a) or (b)) then:
\begin{itemize}
\item[$d$)] $q(S)=0$ unless $p_g(S)=q(S)=1$;
\item[$e$)] unless $K_S$ is ample and $p_g(S)=6,K_S^2=8$ or $p_g(S)=3,K_S^2=2$, 
there is a rational pencil whose general member is a smooth hyperelliptic curve of genus 3 
such that the bicanonical map of $S$ induces the hyperelliptic involution on it.
\end{itemize}
\end{theorem}

We would like to remark that we get the classification
of regular surfaces with $p_g(S)\ge 2$ and non birational bicanonical map. In fact, by the above remark 
and Theorem \ref{mainth} it follows that:
\begin{theorem}
Let $S$ be a smooth minimal surface of general type with $\hbox{q(S)=0}$, $p_g(S)\ge 2$.
Assume that the bicanonical map of $S$ is not birational.

Then if $S$ does not present the standard case it is the smooth minimal model of a Du Val double plane.
\end{theorem}

We remark that Theorem \ref{mainth} also completes the classification of regular surfaces of general type with
$p_g(S)=1$ and non birational bicanonical map.

In fact, in this case if $\varphi_{2K}$ has degree 2 then $S_2$ is a surface of degree $2N-2$ in $\pp\sp N$
and so it is either ruled or a $K3$. The $K3$ case is classified by D.Morrison (\cite{M}).
Otherwise, $\varphi_{2K}$ has degree greater than $2$ and then $K_S^2\le 2$ (cfr. \cite{X2}),
such surfaces are classified by F.Catanese (\cite{C}) for $K_S^2=1$ and by F.Catanese, O.Debarre (\cite{CD})
for $K_S^2=2$.

The paper is organized as follows.
In section 1 we fix some notation and we recall some general facts concerning the surfaces under consideration.
In \S 2 we work out a first easy case, then we prove a result which suffices to get $(b)\Rightarrow(c)$ of
 Theorem \ref{mainth}
and starting from it we prove the implication $(b)\Rightarrow(c)$ in \S 3. 
In \S 4 we prove $(c)\Rightarrow(b)$ and classifying Du Val double planes we get $(d),(e)$.
Finally, in \S 5 we collect some consequences of Theorem \ref{mainth}.

\textbf{Acknowledgement.} I would like to thank prof. Ciro Ciliberto, who suggested the problem, for his
advise and encouragement.
I also would like to
thank Ingrid Bauer and Fabrizio Catanese for fruitful discussions and for their friendly hospitality when I was at
the University of Bayreuth.
I am indebted to Fabrizio Catanese who suggested to remove the hypothesis of regularity in a preliminary version of
the main theorem.

\section{Notation and set up}

Throughout the paper we will mean by surface (resp. curve) a projective algebraic surface (resp. curve)
over the complex
numbers and by a curve on a surface we will mean an effective non zero divisor on the surface.
The symbol $\equiv$ will denote the linear equivalence
 of divisors.

A smooth surface $Y$ is ruled if there exists a surjective morphism $f$ 
onto a curve whose general fibre
is isomorphic to $\pp^1$. If each fibre of $f$ is smooth one says that $Y$ is geometrically ruled.
Let $Y^\prime$ be a singular surface and $Y\rightarrow Y^\prime$ a resolution of the singularities. Then we will say
that $Y^\prime $ is ruled if $Y$ is ruled.

Let $C$ be a reduced curve singular at a point $p\in C$.
The singularity is {\em non essential} if it is:
\begin{itemize}
\item[-] either a double point,
\item[-] or a triple point which resolves to at most a double point after one blow up.
\end{itemize}
otherwise it is {\em essential}.
Let $p^\prime$ be a point infinitely near to $p$.
Then $C$ has an $[r,r]$-point at $(p,p^\prime)$ if it has a point of multiplicity $r$ at $p$ which resolves
to a point of multiplicity $r$ at $p^\prime$ after one blowing up at $p$. We shall denote such singularity
by $[p^\prime\rightarrow p]$.
Notice that an $[r,r]$-point is an essential singularity if and only if $r\ge 3$.

We will use freely the theory of double covers referring to \cite{BPV} for the details.

\subsection{Surfaces with a 2-to-1 rational map.} \label{sect1} Let 
 $S$ be a smooth minimal surface of general type such that there is a
generically finite rational map $\phi:S\dashrightarrow \Sigma$ of degree 2 onto a surface
(for short, a $2$-to-$1$ rational map).

Hence $\phi$ induces an involution $\sigma$ on $S$ which is a morphism since $S$
is minimal of general type. The fixed locus $Fix(\sigma)$ is the union of a smooth reduced curve $R_{\sigma}$ and $k$
distinct points $q_1,..,q_k$.
The canonical projection onto the quotient
$\rho:S\rightarrow \Sigma_{\sigma}:=S/\sigma$ is a double cover, i.e. a finite morphism of degree 2,
branched along the smooth curve $B_{\sigma}=\rho(R_{\sigma})$ and at the points $Q_i=\rho(q_i),\ i=1,\dots ,k$.
The only singularities of $\Sigma_{\sigma}$ are the ordinary double points $Q_1,\dots,Q_k$.

Let $\Hat\pi :\Hat S\rightarrow S$ be the blow-up at $q_1,...,q_k$
and let $E_1,...,E_k$ be the exceptional $(-1)$-curves of $\Hat\pi .$ 
We denote by $\Hat\sigma$ the induced involution
on $\Hat S$ and the quotient $\Hat S/\hat\sigma$ by $\Hat \Sigma.$
Furthermore, we denote $\Hat\pi^{-1}(R_{\sigma})$ by $\Hat R$. Hence $Fix(\Hat\sigma)=\Hat R+E_1+\dots +E_k$
and we get the following commutative diagram
\begin{equation*}
\begin{CD}
\Hat S @>{\Hat\pi }>> S \\
@V{\Hat\rho}VV  @VV{\rho}V \\
\Hat \Sigma @>{\eta }>> \Sigma_\sigma
\end{CD}
\end{equation*}
where the morphism $\eta$ is the minimal resolution of the singularities of $\Sigma_\sigma$ and
$\Hat\rho$ is a double cover branched along the smooth curve $\Hat B=\Hat B^\prime + C_1+\dots +C_k$
where $\Hat B^\prime=\Hat\rho(\Hat R)$ and $C_i=\Hat\rho(E_i), i=1,\dots,k$.
In particular, $C_i=\eta^{-1}(Q_i)$ is a $(-2)$-curve and $\Hat \Sigma$ is smooth.

By the theory of double covers there exists $\Hat \Delta\in Pic(\Hat \Sigma)$ such that $\Hat B\in|2\Hat \Delta|$ and
$\rho\sb \ast \calO\sb {\hat S}=\calO\sb {\hat \Sigma}\oplus \calO\sb {\hat \Sigma}(-\hat\Delta)$. 
Therefore, $K_{\Hat S}=\Hat\rho^\ast(K_{\Hat \Sigma}+\Hat \Delta)$ and we have
\begin{equation*}
\begin{split}
H^i(\Hat S,\calO_{\Hat S}( mK_{\Hat S}))\cong   \\
\cong H^i(\Hat \Sigma, \calO_{\Hat \Sigma}&(m(K_{\Hat \Sigma}+\Hat \Delta)))
\oplus H^i(\Hat \Sigma, \calO_{\Hat \Sigma}(mK_{\Hat \Sigma}+(m-1)\Hat \Delta))
\end{split}
\end{equation*}
for each $i\ge 0$ and $m\ge 0$.

Now we assume that the bicanonical map of $S$ factors through $\phi$, then we have the following commutative diagram
\[
 \xymatrix@!0{
\Hat S \ar [dd]_{\Hat\rho} \ar[rrr]^{\Hat\pi}
   & && S \ar@{.>}[ddrrr]^\phi \ar [dd]_{\rho} \ar@{.>} [rrr]^{\varphi_{2K}} & & & S_2  \\
 \\
\Hat \Sigma \ar[rrr]^{\eta}& && \Sigma_\sigma\ar@{.>} [rrr]^{\eta^\prime}& & & \Sigma \ar@{.>} [uu]_{\phi_2} & &
   }
\]
where $\eta^\prime:=\rho^{-1}\circ \phi$ is a birational map and $\varphi_{2K}$ factors through $\rho$ and $\Hat\rho$.

\begin{remark} In general $\varphi_{2K}$ factors through $\rho$ if and only if either
\newline $H^0(\Hat \Sigma,\calO_{\Hat \Sigma}(2K_{\Hat \Sigma}+2\Hat \Delta))=0$ or
$H^0(\Hat \Sigma,\calO_{\Hat \Sigma}(2K_{\Hat \Sigma}+\Hat \Delta))=0.$
\end{remark}
Therefore, we know that in our situation one of the above vector spaces has to be trivial. 
In fact, in the following refined version of a proposition by M.Mendes Lopes and R.Pardini
(cfr.\cite{MP2}, Proposition 2.1) 
we will see that in our situation $H^0(\Hat \Sigma,\calO_{\Hat \Sigma}(2K_{\Hat \Sigma}+\Hat \Delta))=0.$

\begin{proposition}\label{wellknow}
Let $S$ be a smooth minimal surface of general type and $\sigma$ an involution acting on $S$.
Let $\hat S$ be the blow up of $S$ at the isolated fixed point of $\sigma$ and
$\hat \rho:\hat S\rightarrow \hat\Sigma:=\hat S/\hat \sigma$ the canonical projection onto the quotient.
Denote by $\hat \Delta\in Pic(\hat\Sigma)$ a divisor such that
$K_{\hat S}=\hat \rho^\ast(K_{\hat\Sigma}+\hat\Delta)$ and by $k$ the number of isolated fixed points of $\sigma$.
Then 
\begin{itemize}
\item [a)] $h^i(\Hat \Sigma,\calO_{\Hat \Sigma}( 2K_{\Hat \Sigma}+ \Hat \Delta))=0,\ for\ i>0$;
\item [b)] let $R_{\sigma}$ be the divisorial part of $Fix(\sigma)$, then
\begin{equation*}
\begin{split}
i)\ \ \ k&=K_S^2-2\chi(\calO_S)+6\chi(\calO_{\hat\Sigma})-2h^0(\Hat \Sigma,\calO_{\Hat \Sigma}( 2K_{\Hat \Sigma}+ \Hat \Delta)\\
ii)\ \ \ k&=K_S.R_{\sigma}-4\chi(\calO_S)+8\chi(\calO_{\hat\Sigma})
\end{split}
\end{equation*}
\item [c)] Assume that $p_g(\hat\Sigma)=0$, then the following three conditions are equivalent
\begin{itemize}
\item[$i$)] the bicanonical map of $S$ factors through $\hat\rho$;
\item[$ii$)] $h^0(\Hat \Sigma, \calO_{\Hat \Sigma}(2K_{\Hat \Sigma}+ \Hat \Delta))=0$;
\item[$iii$)] $k=K_S^2-2\chi(\calO_S)+6\chi(\calO_{\hat\Sigma})$.
\end{itemize}
\end{itemize}
\end{proposition}
\begin{proof}
$a)$ We use the notation introduced before. As $\Sigma_\sigma$ has at most canonical singularities, we have that
 $2K_S=\rho ^*(2K_{\Sigma_\sigma}+R_\sigma)$.
Therefore, $2K_{\Sigma_\sigma}+R_\sigma$ is nef and big because
$2K_S$ is nef and big, and so $2K_{\Hat \Sigma} + \Hat B=\eta \sp *(2K_{\Sigma_\sigma}+R_\sigma)$ is nef and big.

 On the other hand we have the following equality of $\mathbb Q $-divisors
\[
K_{\Hat \Sigma}+\Hat \Delta=\frac12 (2K_{\Hat \Sigma} +\Hat B)+\frac12 \sum C_j
\]
where $\frac12 \sum C_j$ is an effective $\mathbb Q$-divisor with zero integral part.
Hence by the Kawamata-Viehweg vanishing theorem it follows that  $h\sp i(\Hat \Sigma,2K_{\Hat \Sigma}+\Hat \Delta)=0$, $i>0$.

$b)$ By $a)$ and the Riemann-Roch formula we get:
\[
h^0(\Sigma,\calO_{\Hat \Sigma}(2K_{\Hat \Sigma}+\Hat\Delta))=\chi(2K_{\Hat \Sigma}+\Hat \Delta)=
\chi(\calO_{\hat\Sigma})+K_{\Hat \Sigma}^2+\frac32 K_{\Hat \Sigma}.\Hat \Delta+\frac12\Hat \Delta^2
\]
and
\[
\chi(\calO_{\Hat S})=\chi (\calO_{\Hat \Sigma})+\chi(K_{\Hat \Sigma}+\Hat \Delta)=
2\chi(\calO_{\hat\Sigma})+\frac12(\Hat \Delta^2+\Hat \Delta.K_{\Hat \Sigma})
\]
since  $\rho_\ast\calO_{\hat S}=\calO_{\hat\Sigma}\oplus \calO_{\hat\Sigma}(-\hat\Delta)$.
On the other hand we have
\[
\left\{\begin{array}{l}
k=K_S\sp 2-K_{\Hat S}^2=K_S^2-2(K_{\Hat \Sigma}+\Hat \Delta)^2 \\
R_\sigma.K_S-k=(\Hat R+\sum_{i=1}^k E_i).K_{\Hat S}=2\Hat \Delta.(K_{\Hat \Sigma}+\Hat \Delta)
\end{array}
\right.
\]
and so using the above equalities we get
\[
\left\{\begin{array}{l}
k=K^2_S-2\chi(\calO_S)+6\chi(\calO_{\hat\Sigma})-2h^0(\hat \Sigma, \calO_\Sigma(2K_{\Hat \Sigma}+\Hat \Delta)) \\
k=R_\sigma.K_S-4\chi(\calO_S)+8\chi(\calO_{\hat\Sigma})
\end{array}
\right.
\]

$c)$ First of all recall that we have
\[
p_g(S)=p_g(\Hat S)=h^0(\Hat \Sigma,\calO_{\Hat \Sigma}(K_{\Hat \Sigma}+ \Hat \Delta))+h^0(\Hat \Sigma,\calO_{\Hat \Sigma}(K_{\Hat \Sigma}))
=h^0(\Hat \Sigma,\calO_{\Hat \Sigma}(K_{\Hat \Sigma}+\Hat \Sigma)).
\]

Therefore, if $p_g(S)>0$ there is a non zero effective divisor $2D\in|2K_{\Hat \Sigma}+\Hat B|$
where $D\in |K_{\Hat \Sigma}+\Hat \Delta|$.
Whence, if $p_g(S)>0$
the bicanonical map of $S$ factors through $\hat \rho$ if and only
if $h^0(\Hat \Sigma,\calO_{\Hat \Sigma}( 2K_{\Hat \Sigma}+\Hat \Delta))=0$.

If $p_g(S)=0$ then from $b)$ it follows 
\[
k=R_\sigma.K_S+4=K^2_S+4-2h^0(\hat \Sigma, \calO_\Sigma(2K_{\Hat \Sigma}+\Hat \Delta)) \\
\]
where $k\ge 4$ since $K_S$ is $nef$.
Now assume that $\varphi_{2K}$ factors through $\hat \rho$ and that
$h^0(\Hat \Sigma,\calO_{\Hat \Sigma}(2K_{\Hat \Sigma}+\Hat \Delta))\neq 0$.
Then $h^0(\Hat \Sigma,\calO_{\Hat \Sigma}(2K_{\Hat \Sigma}+\Hat B))=0$ and we have
\[
h^0(\Hat \Sigma,\calO_{\Hat \Sigma}(2K_{\Hat \Sigma}+\Hat \Delta))=h^0(\Hat S,\calO_{\Hat S}(2K_{\hat S}))=h^0(S,\calO_S(2K_{S}))
= K_S^2+1
\]
which by the above equality implies
\[
k=K_S^2+4-2(K_S^2+1)=-K_S^2+2\le 1.
\]
A contradiction.
Whence, the bicanonical map of $S$ factors through $\hat\rho$ if and only if
$h^0(\Hat Y,\calO_{\Hat \Sigma}( 2K_{\Hat \Sigma}+\Hat \Delta))=0$ and, by $b,i)$,
 the equivalence with $c,iii)$ is clear.
\end{proof}
\subsection {Canonical resolution.} (cfr. \cite{H}, \cite{BPV})
Let $W_0$ be a smooth surface. Assume that there exists a double cover $S_0 \rightarrow W_0$
branched along a reduced curve
$B_0\subset W_0$. Then $S_0$ is normal and it is smooth if and only if $B_0$ is smooth. If $S_0$ is singular the
singularities of $S_0$ can be resolved in a natural way by the $canonical$ $resolution$. Briefly, there is a
commutative diagram

\[
 \xymatrix@!0{
S^\ast \ar [dd]_{\rho^\prime} \ar[rrrrrr] & &&&& & S_0 \ar [dd]  \\
 \\
W_s \ar[rr]^{\omega_{s}}& & W_{s-1}\ar[rr]^{\omega_{s-1}} & &\dots \ar[rr]^{\omega_{1}}& & W_0
   }
\]
such that, for each $i=0,\dots,s-1,$ $\omega_{i+1}$ is the blow up of $y_i\in W_i$ , where $y_i\in B_i$ is a singular point
 of $B_i$.
 
Let $m_i$ be the multiplicity of $B_i$ at $y_i$ and $\calE_{i+1}=\omega_{i+1}^{-1}(y_i)$
the exceptional curve of $\omega_{i+1}$, hence $B_{i+1}=(\omega_{i+1})^*(B_i)-2[\frac{m_i}2]\calE_{i+1}$ where
$[\frac{m_i}2]$ is the greatest integer lesser than or equal to
$\frac{m_i}2$.
Furthermore, the curve $B_s\subset W_s$ is smooth, $\rho^\prime$ is a double cover branched along
$B_s$  and $S^\ast \rightarrow S_0$ is a birational morphism.

Let us denote by $\omega=\omega_1\circ\dots\circ\omega_s$ the composition and by
$\calE_{i+1}^\ast=\omega^\ast(y_i)$
the exceptional $(-1)$-cycle with reduced support $\omega^{-1}(y_i)$.
Hence the following equalities hold
\begin{align*}
K_{W_s}=\omega^\ast(K_{W_0})+\sum\calE_i^\ast; \ \ \ \ \ \ &\ \ \ \ \ \ \ \calE_j^\ast.\calE_h^\ast=-\delta_{j,h}; \\
B_s=\omega^\ast(B_0)-\sum &2\left[\frac{m_i}2\right]\calE_i^\ast 
\end{align*}
where $\delta_{j,h}$ is the Kronecker symbol.

Notice that $S^\ast$ is also the canonical resolution of the double cover $S_i\rightarrow W_i$ branched along $B_i$, for
each $i=1,\dots,s$.

\begin{lemma}\label{lemmaFactor}
Let $S,\hat S$ and $\hat \Sigma$ be as in Proposition \ref{wellknow}.
Let $\psi:\Hat \Sigma\rightarrow \Sigma^\prime$ be a birational morphism onto a smooth surface and
consider a factorization of $\psi$ in blow ups
\[
\psi :\Hat \Sigma=\Hat \Sigma_0 \stackrel{\psi_1}\longrightarrow \Hat \Sigma_1 \stackrel{\psi_2}\longrightarrow \dots
\stackrel{\psi_t}\longrightarrow \Hat \Sigma_t=\Sigma^\prime
\]
For $i=1,\dots ,t$, denote by $ \Hat y_i\in\hat\Sigma_i$ the center of the blow up $\psi_i$ and by
 $\Hat \calE_i=\psi_i^{-1}(\hat y_i)$
 the exceptional curve of $\psi_i$. Moreover, let 
$\hat B_i$ be the image of $\Hat B $ in $\Hat \Sigma_i$. Set $\Hat B_0=\Hat B$.
Then for each $i\ge 1$
\begin{enumerate}
\item if $\Hat B_i$ has a singularity at a point $z$ then either $\Hat y_i=z$ or there exists 
$j < i$ such that $(\psi_i\circ \dots \circ \psi_{j+1} )(\Hat y_j)=z$;
\item $\Hat \calE_i$ belongs to $\Hat B_{i-1}$ if and only if the multiplicity of $\Hat B_i$ at $\Hat y_i$ is odd;

\item $\Hat B_i$ is singular at $\Hat y_i$;

\item $\Hat S$ is the canonical resolution of the double cover of $\Sigma^\prime$ branched \hbox{along $\hat B_t$}. 
\end{enumerate}\end{lemma}
\begin{proof} We keep the notation from section $1.1$.
Since $\hat B$ is smooth, $1$) is clear.
2) For $i=1,\dots,t,$ let $\Hat \Delta_i\in Pic(\hat \Sigma_i)$ denote $(\psi_i\circ \dots \circ \psi_1)_\ast(\Hat \Delta)$. 
If $\Hat \calE_i\not \subset\Hat B_{i-1}$ 
the multiplicity of $\Hat B_i={\psi_i}_\ast(\Hat B_{i-1})$ at $\Hat y_i$ is
$\Hat \calE_i.\Hat B_{i-1}=2\Hat \calE_i.\Hat \Delta_{i-1}$, an even number.
On the other hand if $\Hat \calE_i\subset \Hat B_{i-1}$ we have
$\Hat\calE_i.(\Hat B_{i-1}-\Hat \calE_i)=2\Hat \calE_i.\Hat \Delta_{i-1}-1$, and so the
multiplicity of $\Hat B_i={\psi_i}_\ast(\Hat B_{i-1})={\psi_i}_\ast(\Hat B_{i-1}-\Hat\calE_i)$ at $\Hat y_i$ is odd.

3) Let $\calE\subset \hat \Sigma$ be a $(-1)$-curve and $E\subset \hat S$ a reduced and irreducible
curve such that $\hat \rho(E)=\calE$.
If $\calE\subset \hat B$ then $E^2=\frac12\calE^2=-\frac12$, a contradiction.
If $\calE\cap\hat B=\emptyset$ then $E^2=-1$ and $E.E_i=0,i=1,\dots,k,$ hence $\hat\pi(E)\subset S$ is a $(-1)$-curve, a
contradiction. Therefore, $\calE\not\subset \hat B$ and $\calE.\hat B\ge 1$, that is $\calE.\hat B\ge 2$ as
$\hat B\equiv 2\hat \Delta$. 
In particular, it follows that $\hat B_1$ is singular at $y_1$.

Now assume $i>1$. By $1),2)$ and the inductive hypothesis, $\hat \calE_i\not \subset \hat B_{i-1}$ implies that $\hat B_i$
has multiplicity $\hat \calE_i.\hat B_{i-1}\ge 2$ at $\hat y_i$ while for $\hat \calE_i\subset \hat B_{i-1}$
we get \hbox{$\hat \calE_i\cap(\hat B_{i-1}-\hat \calE_i)\neq \emptyset$}.

Hence $\hat B_i$ is singular at $\hat y_i$ if $\hat \calE_i\not \subset \hat B_{i-1}$ and
$\hat B_i={\psi_i}_\ast(\hat B_{i-1}-\hat \calE_i)$ has multiplicity $\ge 1$ at $\hat y_i$ if
$\hat \calE_i\subset \hat B_{i-1}$.
In the second case if $\hat \calE_i.(\hat B_{i-1}-\hat \calE_i)=1$ we can assume 
$\{\hat y_{i-1}\}=\hat \calE_i\cap(\hat B_{i-1}-\hat \calE_i)$ and still by induction we get that the strict transform
$\st \calE_{i-1}$ (resp. $\st \calE_i$) of $\hat \calE_{i-1}$ (resp. $\hat \calE_i$) on $\hat\Sigma$
is a $(-1)$-curve ($(-2)$-curve)
belonging (do not belonging) to $\hat B$ such that $\st \calE_i.\st \calE_{i-1}=1$ and $\hat B.\st \calE\sb{i-1}=2$.
Therefore, taking the pull back
to $\hat S$ of $\st \calE_{i-1}$ and then pushing it down to $S$ we get a smooth rational curve with 
selfintersection greater than or equal to $-1$. A contradiction.

Finally, for $4)$ it is easily seen that, since $\hat B$ is smooth, $1),2),3)$ characterize the canonical resolution of
the double cover of $\Sigma_t$ branched \hbox{along $\Hat B_t$.}  
\end{proof}

\section{Proof of Theorem \ref{mainth}: part I}

In this section and in the next one we will prove the implications  $(a) \Rightarrow (b)$, $(a) \Rightarrow (c)$ of
Theorem \ref{mainth}.
Hence, throughout these two sections we will assume that $S$ is a smooth minimal surface of general type  such that 
the bicanonical map factors through a $2$-to-$1$ map $\phi:S\dashrightarrow \Sigma$ 
onto a rational or ruled surface.
We also assume that $S$ does not present the standard case, in particular $K_S^2\le 9$.

Therefore, from section \ref{sect1} we get the commutative diagram
\[
 \xymatrix@!0{
 \Hat S\ar [dd]_{\Hat\rho} \ar[rrr]^{\pi} & && S  \ar[dd]_{\rho} \ar[rrr]^{\varphi_{2K}} \ar@{.>}[rddrr]^{\phi}& && S_2  \\
 \\
\Hat \Sigma \ar[rrr]^{\Hat\eta} & && \Sigma_\sigma \ar@{.>}[rrr]^{\eta} & && \Sigma\ar@{.>}[uu]_{\phi_2}
}
\]
where $\Hat \Sigma$ is a rational or ruled surface since $\eta$, $\eta^\prime$ are birational maps.
In particular, as $\Hat \Sigma$ is smooth it is either ruled or $\pp^2$.
\begin{proposition}
If $\Hat\Sigma \cong \pp^2$  then $q(S)=0$ and $\Hat B$ is a smooth curve of degree 8 or 10.
We have respectively $p_g(S)=3, K_S^2=2$ and $p_g(S)=6$, $K_S^2=8$.
Moreover, $K_S$ is ample.
\end{proposition}

\begin{proof} First of all notice that the involution $\sigma$ induced by $\phi$  on S does not have isolated
fixed points, otherwise there would be some $(-2)$-curve contained in $\pp^2$ (cfr. (1.1)).

Hence $\Hat S=S$, $\Hat \Sigma=\Sigma_\sigma$ and $\rho $ is a (finite)
double cover. Therefore, $\Hat B$ is smooth and denoting by $2d$ the degree of $\Hat B\equiv 2\Hat \Delta$ we get:
\begin{align*}
9&\geq K_S^2=2(K_{\pp^2}-\Hat \Delta)^2=2(d-3)^2, \\
h&^i(S,\calO_S(K_S))=h^i(\pp^2,\calO_{\pp^2}(d-3))+h^i(\pp^2,\calO_{\pp^2}(-3)),
\end{align*}
hence $2d\leq 10$. On the other hand we have $2d\geq 8$, since $S$ is of general type.
So $3\le d\le 4$ and $q(S)=0,\ p_g(S)=\frac12 d(d-3)+1$.

We notice that there cannot be a $(-2)$-curve on $S$, since it would map to a $(-1)$-curve or a $(-2)$-curve in $\pp^2$,
whence $K_S$ is ample.     \end{proof}

From now on we will assume that $\Hat \Sigma$ is ruled. Let $\Sigma_e$ be a geometrically ruled surface. 
We denote by $C_0$ a section of $\Sigma_e$ such that the self intersection
$C_0^2=-e\le 0$ is the smallest possible and by $\Gamma\cong \pp^1$ we denote a fibre of the ruling.
Recall that $C_0$ and $\Gamma $ generate $Pic(\Sigma_e)$.

Hence there is a birational morphism $\varphi:\hat \Sigma\rightarrow \Sigma_e$ and 
setting $B=\varphi_\ast(\Hat B)$ we can write
\[
B\equiv\xi C_0 + (\frac12 \xi e +\zeta)\Gamma
\]

Following Xiao \cite{X2}, we can assume $\varphi$ to be such that
\begin{itemize}
\item[$\dagger$)] $\xi= B.\Gamma$ is minimal;

\item[$\ddagger)$] the greatest multiplicity of the singularities of $B$ is minimal, and the number of singularities of
$B$ with the greatest multiplicity is minimal, among all the choices satisfying condition $(\dagger)$;
\end{itemize}
where an [r,r]-point is considered as a unique singularity of multiplicity strictly between $r$ and $r+1$.

\textbf{Remark.}  Let $\Hat H$ be the pull back to $\hat S$ of a general $\Gamma\in|\Gamma|$.
Hence $\varphi\circ\Hat\rho|_{\Hat F}:\Hat F\rightarrow  \Gamma$ is a double cover branched in $\Gamma.B$ points.
Therefore, $|\Hat H|$ is a pencil of curves of genus $\frac12(\Gamma.B-2)$.
In particular, we assume $\xi\geq 8$ since $S$ does not present the standard case.

The main result of this section is the following:
\begin{theorem}\label{main2}
Let $S$ be a smooth minimal surface of general type does not presenting the standard case
and $\sigma$ an involution acting on $S$ such that the quotient $S/\sigma$ is a ruled surface.
Let $\hat S$ be the blow up of $S$ at the isolated fixed point of $\sigma$ and
$\hat \rho:\hat S\rightarrow \hat\Sigma=\hat S/\hat \sigma$ the projection onto the quotient.
Let $\varphi :\Hat \Sigma \rightarrow \Sigma_e$ be a birational morphism having the
properties $\dagger)$ and $\ddagger)$.

Assume that the bicanonical map of $S$ factors through $\hat \rho$.
Then $\Sigma_e$ is rational and only the following possibilities can occur:
\begin{itemize}
\item [i)] $\xi =8,\ \zeta=6$;
\item [ii)] $\xi = 8,\ \zeta=8+2i$, where $1\leq i\leq 5$. The essential singularities of $B$ are: $i+1$
     $[5,5]$-points, possibly some $4$-tuple points or $[3,3]$-points.
\end{itemize}
\end{theorem}

\begin{remark}\label{xlist} The idea of this theorem goes back to Xiao Gang. In fact, in \cite{X2} Proposition 6
he proves a weakly result,
namely:\newline $a)$ he further assumes the bicanonical map to be $2$-to-$1$ onto a ruled surface and 
that $h^0(\hat\Sigma, \calO_{\Hat \Sigma}(2K_{\Hat\Sigma}+\Hat \Delta))=0$;
\newline $b)$ he claims that under these hypotheses 
$\Sigma_e$ is rational and only the following possibilities can occur: $(i), (ii)$ as above and 
\begin{itemize}
\item [$iii$)]  $\xi = 12,\ \zeta=14$, and $B$ has three [$7,7$]-points, possibly some non essential singularities;
\item [$iv$)]    $\xi = 16,\ \zeta=18$, and $B$ has three [$9,9$]-points, an $8$-tuple point,
 possibly some non essential singularities.
\end{itemize}
\end{remark}
\begin{remark} In fact, Theorem \ref{main2} suffices to prove implication $(a)\Rightarrow (b)$ of
\hbox{Theorem \ref{mainth}}. In particular, we have that $\Sigma_e$ is the Hirzebruch surface 
$\ff_e=\textbf{P}(\calO_{\pp^1}\oplus\calO_{\pp^1}(-e))$. 

\end{remark}
\begin{remark}\label{mainremark} We will prove the above theorem in several steps:\newline
1) we remark that looking carefully at the Xiao's proof it is easy to see that the argument still works 
if one suppose that the bicanonical map factors through a rational map of
degree two onto a rational or ruled surface;\newline
2) moreover, in our situation  we have that
$h^0(\hat\Sigma, \calO_{\Hat \Sigma}(2K_{\Hat\Sigma}+\Hat \Delta))=0$ by \hbox{Proposition \ref{wellknow}};
\newline
3) therefore, we are now reduced to prove the following proposition:
\end{remark}
\begin{proposition}\label{not}
In the hypotheses of Theorem \ref{main2}, cases $(iii)$, $(iv)$ above do not occur.
\end{proposition}

{\bfseries Notation.} From now on we will refer to $\Hat S$ as a surface of type $\mathcal{S}_I$
(resp. $\mathcal{S}_{II}$, $\mathcal{S}_{III}$, $\mathcal{S}_{IV}$)
meaning that we consider $\Hat S$ associated to the commutative diagram
\[
\xymatrix@!0{
&\Hat S\ar [dd]_{\Hat\rho} \ar[rrrdd]   \\
 \\
\Hat B\subset&\Hat \Sigma \ar[rrr]^{\varphi} & && \ff_e & \supset B
}
\]
such that the morphism $\varphi :\Hat\Sigma \rightarrow \ff_e$ has the properties $(\dagger)$, $(\ddagger)$
and $B=\varphi_\ast(\Hat B)$ is as in Proposition \ref{main2},(i) (resp. \ref{main2},(ii),
Remark \ref{xlist},(iii),(iv)).

Let $p\in\ff_e$ be a point. We denote by $elm_p$ the elementary transformation centered at $p$,
that is the result of blowing up $p$ and then contracting the fibre of
the ruling passing through $p$.

\begin{lemma}\label{corollMain}
Let $\Hat S$ be a surface of type $\mathcal{S}_{II}$ (resp. $\mathcal{S}_{III}$, $\mathcal{S}_{IV}$).
Let $[p^\prime\rightarrow p]$ be an $[r,r]$-points of $B$ such that $r=5$ (resp. $7,9$).
Let $\Gamma_p\in |\Gamma| $ be the fibre such that $p\in \Gamma_p$.

Then, $p^\prime$ is infinitely near to $\Gamma_p$, $\Gamma_p$ belongs to $B$ and
two distinct singular $[r,r]$-points lie on distinct fibres. Finally, $C_0$ does not belong to $B$.
\end{lemma}
\begin{proof}
Suppose that $p^\prime$  is not infinitely near to $\Gamma_p$.
Then $\Gamma_p\not \subset B$ and since $r=\frac12\xi+1=\frac12B.\Gamma_p+1$ we can apply $elm_{p^\prime} \circ elm_p$ to obtain a new model with less
singularities of maximal multiplicity. A contradiction.

Now the other claims are clear. 
\end{proof}

\begin{lemma}\label{lemmaE}
Let $\Hat S$ be a surface of type $\mathcal{S}_{II}$ or $\mathcal{S}_{III}$ or $\mathcal{S}_{IV}$.
Then
\begin{itemize}
\item[-] if $\Hat S$ is of type $\mathcal{S}_{II}$ then $0\leq e \leq \frac {7+i}{4}, \ 1\leq i \leq 5$;
\item[-] if $\Hat S$ is of type $\mathcal{S}_{III}$ or $\mathcal{S}_{IV}$ then we can assume $e=1$.
\end{itemize}
\end{lemma}
\begin{proof}
By Lemma \ref{corollMain} the curve $B$ contains $i+1$ (resp. $3$) fibres if $\Hat S$ is of type $\mathcal{S}_{II}$
(resp. $\mathcal{S}_{III}$ or $\mathcal{S}_{IV}$) and $C_0\not \subset B$.

Therefore, if $\Hat S$ is of type $\mathcal{S}_{II}$ we get
\[
i+1\leq C_0.B=-4e+8+2i
\]
that is
\[
e\leq \frac {7+i}{4}
\]
and analogously we get $e< 2$ if $\Hat S$ is of type $\mathcal{S}_{III}$ or $\mathcal{S}_{IV}$.

Let us now suppose that $\Hat S$ is of type $\mathcal{S}_{III}$ (resp. $\mathcal{S}_{IV}$)
and $e=0$. Then we can choose $C_0 $ such that there exists
a [7,7]-point (resp. [9,9]-point), say [$p^\prime \rightarrow p$],
such that $p\in C_0$. Now performing $elm_p$ we get a model with $e=1$
and the same singularities.
\end{proof}

For the remainder of the section we will assume  that $\Hat S$ is of type $\mathcal{S}_{III}$ or $\mathcal{S}_{IV}$.
Therefore, by Lemma \ref{lemmaFactor} $\Hat S$ is the canonical resolution of the double cover
of $\ff_1$ branched along
a reduced curve $B=B^\prime+\Gamma_1+\Gamma_2+\Gamma_3$ such that $\Gamma_i\in|\Gamma|$ and
$B^\prime\in |12C_0+17\Gamma|$ (resp. $B^\prime\in |16 C_0+21\Gamma|)$.

We will denote by $[p^\prime_i \rightarrow p_i]$ the $[7,7]$-point (resp. $[9,9]$-point)
of $B$ such that $p_i\in L_i,\ i=1,2,3$.

\begin{lemma}
$(a)$ For any curve $C\subset \ff_1$ sitting in the linear system $|C_0+\Gamma|$
we have that $\{p_1,p_2,p_3\}\not \subset C$;

$(b)$ $p_i\not \in C_0,\ i=1,2,3$.
\end{lemma}
\begin{proof}
We give the proof for the case $B^\prime\in |12C_0+17\Gamma|$, the other one is completely analogous.
Suppose that there exists $C\in|C_0+\Gamma|$ such that $\{p_1,p_2,p_3\}\subset C$. Then
$21\le C.B=(C_0+\Gamma).(12C_0+20\Gamma)=20$ implies that $C$ belongs to $B$.
Hence $C$ is tangent to $\Gamma_i$ at $p_i$,
for $i=1,2,3,$ a contradiction.
Analogously we get $p_i\not\in C_0$.
\end{proof}

Since the $[p_i^\prime\rightarrow p_i]$'s are the singularities of $B$ with maximal multiplicity
we have a factorization $\varphi=\pi_1\circ\pi_2\circ\dots \circ\pi_6\circ \varphi^\prime $
such that $\pi_i$ (resp. $\pi_{3+i}$) is the blow up at $p_i$ (resp. $p_i^\prime$), $i=1,2,3$.
We set $\Hat W=\varphi^\prime (\Hat \Sigma)$.

Let $\pi_0:\ff_1\rightarrow \pp^2$ be the morphism contracting $C_0$ to a point $p_0\in\pp^2$.
Then $L_i=\pi_0(\Gamma_i)$ is a line passing through $p_0$ and by the above lemma
$p_1,p_2,p_3,$ are non collinear points such that  $p_i\neq p_0,\ i=1,2,3.$ (by abuse of notation we denote by the same
letter the image of $p_i$ in $\pp^2$). Hence we have the commutative diagram
\[
\xymatrix@!0{
\Hat \Sigma \ar[ddrrrr]^\varphi \ar[dd]_{\varphi^\prime}     \\
\\
\hat W \ar[ddrrrr]^\pi \ar[rrrr]^{\pi_1\circ\dots\circ\pi_6\ \ \ } &    &   &    & \ff_1 \ar[dd]^{\pi_0} &   \\
\\
& & & & \pp^2 &
}
\]
where $\pi=\pi_0\circ\pi_1\circ\dots \circ\pi_6:\Hat W\rightarrow \pp^2$ is the composition.
We set $\calE_0^\ast=\pi^\ast(p_0)$, $\calE_i^\ast=\pi^\ast(p_i)$, $\calE_{3+i}^\ast=\pi^\ast(p_i^ \prime)$, $i=1,2,3$.

\begin{lemma}\label{lemma}
Let $L$ be a line in $\pp^2$.
Then
\[
h^0(\Hat W,\calO_{\Hat W}(\pi^\ast (5L) -\calE_0^\ast-\sum_{1=1}^6 2\calE_i^\ast))=2
\]
and the general element
\[
D\in |\pi^*(5L) -\calE_0^\ast-\sum_{1=1}^6 2\calE_i^\ast|
\]
is a smooth and irreducible rational curve on $\Hat W$ such that $D^2=0$.

In particular, $|D|$ defines surjective morphism $f:\Hat W\rightarrow \pp^1$ such that the general fibre is
isomorphic to $\pp^1$.
\end{lemma}
\begin {proof}
Let $\mathcal{C}_i$ be the conic in $\pp^2$ passing through $p_1,p_2,p_3,$ tangent to $L_j,L_k$, $\{i,j,k\}=\{1,2,3\}$.
 Then $\mathcal{C}_i$ is smooth since $p_1,p_2,p_3$ are not collinear and $p_i\neq p_0$.
Let  $\st D_i\subset \Hat W$ be the strict transform of the curve
$D_i=2\mathcal{C}_i+L_i,\  i=1,2,3$.
Hence $\st D_i,\st D_j$ do not have common components if
$i\ne j$ and $\st D_i\in |D|$.
In particular, $h^0(\Hat W,\calO_{\Hat W}(D))\ge 2$ and  $|D|$ does not have fixed part.

A straightforward calculation yields $D^2=0,D. K_{\Hat W}=-2.$
Therefore, the rational map $f$ defined by $|D|$ is a surjective morphism onto a curve and
$D\in |a\mathcal L|$, where $\mathcal L$ is a general fibre of $f$.

On the other hand $\calE_0\sp\ast.D=\calE_0^\ast.(\pi^*(5L) -\calE_0^\ast-\sum_{1=1}^6 2\calE_i^\ast)=1$. 
Therefore, $a=1$ and a standard argument completes the proof.  
\end{proof}

Now we are ready to prove Proposition \ref{not} and then Theorem \ref{main2} by \hbox{Remark \ref{mainremark}}.\newline
\begin{proof}[Proof of Proposition \ref{not}.]
Let $D$ be as in the above lemma.
Denote by $\pi^\prime=\pi_1\circ\dots\pi_6:\hat W\rightarrow \ff_1$ the composition 
such that $\varphi=\pi^\prime \circ \varphi^\prime$. 
 Hence, a straightforward calculation yields $D\in |{\pi^\prime}^\ast(4C_0+5\Gamma)-\sum_{1=1}^6 2\calE_i^\ast|$.

 Assume that $\Hat S$ is of type $\mathcal S_{III}$. Then $B\equiv 12C_0+20\Gamma$ and $\Hat S$ is 
the canonical resolution of the double cover of $\ff_1$ branched along $B$. 
Hence we have
\[
\varphi^\prime_\ast(\Hat B)={\pi^\prime}^\ast(12C_0+20\Gamma)-\sum_{1=1}^3 6\calE_i^\ast-\sum_{1=4}^6 8\calE_i^\ast
\]
since the $[p_i^\prime\rightarrow p_i]$'s are $[7,7]$-points.
Therefore, we get
\[
D.(\varphi_\ast(\Hat B))=(4C_0+5\Gamma).(12C_0+20\Gamma)-12\cdot 3-16\cdot 3=8
\]
which is a contradiction. Indeed in this case
there is a birational morphism $\tilde \varphi:\Hat W\rightarrow \ff_e$ such
that $\tilde \varphi_\ast(D)$ is a ruling of $\ff_e$ and so if we consider the morphism
$\tilde \varphi\circ\varphi^\prime:\Hat \Sigma\rightarrow \ff_e$ we get
$((\tilde \varphi\circ\varphi^\prime)_\ast(\Hat B)).(\tilde \varphi_\ast(D))=8<12$. But we are assuming
that $\varphi$ has the property $(\dagger)$.

An analogous argument shows that $\Hat S$ can not be of type $\mathcal S_{IV}$. 
\end{proof}

\section{Proof of Theorem \ref{mainth}: part II}\label{sectDV}
In this section we will prove that if $\hat S$ is of type $\mathcal S_I$ (resp. $\mathcal S_{II}$) then $S$
is the minimal model of a Du Val double plane. Therefore, we get the implication $(a)\Rightarrow (c)$
of Theorem \ref{mainth}.

\begin{proposition}
Let $S$ be such that $\Hat S$ is a surface of type $\mathcal{S}_I$. Then there exists a birational morphism
$\psi:\hat\Sigma\rightarrow X$ onto as smooth surface
such that setting $G:=\psi\sb\ast(\hat B)$ we have:
\begin{itemize}
\item[$a$)] either $X\cong \pp^2 $ and $G$ is a reduced curve of degree $10$ with possibly some $[3,3]$-points
and no other essential singularities;

\item[$b$)] or $X\cong \ff_2$ and $G=C_0+G'$, where $G^\prime\cap C_0=\emptyset$ and
$G'$ is a reduced curve in the linear system $|7C_0+14\Gamma|$ with
at most non essential singularities;

\item[$c$)]or $X\cong \pp^2 $ and $G=G^\prime+L_1$ where $L_1$ is a line and $G^\prime $ is a reduced curve of
 degree $11$.
 In this case
$G$ has the following essential singularities: a $4$-tuple point and a $[5,5]$-point on $L_1$, possibly some
$[3,3]$-points.
\end{itemize}

  In particular, $\hat S$ and $S$ are respectively the canonical resolution and the smooth 
minimal model of a Du Val double plane.
\end{proposition}
\begin{proof}
This was already partially proved by Xiao Gang. In fact, we have a morphism $\varphi:\hat\Sigma \rightarrow \ff_e$
 such that
$B\equiv 8C_0+(4e+6)\Gamma$ and by \cite{X2} \hbox{Proposition 7}, either $e=1$ or $e=2$, hence either 
$B\equiv 8C_0+10\Gamma$ or
$B\equiv 8C_0+14\Gamma$. Still by \cite{X2} Proposition 7 the essential singularities of $B$  are possibly
4-tuple points or [3,3]-points.

If $e=1$ let $cont_{C_0}$ be the morphism which contracts the $(-1)$-section to a point $p\in\pp^2$. 
We obtain a morphism onto $\pp^2$
\[
\psi:=cont_{C_0}\circ \varphi:\hat\Sigma \rightarrow \pp^2
\]
such that $G:=\psi\sb\ast(\hat B)$ is a curve of degree $10$. Notice that $C_0.B=2$ and so $G$ has
either a triple point or a double point at $p$ depending on
$C_0\subset B$ or not. Thus the essential singularities of $G$ are possibly 4-tuple points or
[3,3]-points.
 Suppose that $G$ has a $4$-tuple point at say $q$. Then the pull-back to $\hat S$
of the pencil $|L-q|$ of lines through $q$ is a pencil of curves of genus 2, a contradiction.

If $e=2$ we have $C_0.B=-2$, thus $C_0\subset B$ and
$B=C_0+ B'$ where $B':=B-C_0$ is a curve such that
$B'\cap C_0=\emptyset$ since $B$ is reduced.

Suppose that $B$ has two 4-tuple points, say $p,q$. Then the pull-back to $\hat S$ of the pencil
$|C-p-q|$, where $C\equiv C_0+2\Gamma$, is a pencil of curves of genus two, a contradiction.

If $B$ has only non essential singularities we set $X=\ff_2$ and $G=B$ (case (b)).
If $B$ has a $4$-tuple point at say $p$ we consider the projection from $p\in\ff_2$ onto the plane, 
i.e. perform an elementary transformation centered at $p$ and then contract the proper transform of $C_0$.
Since we blow up a singular point of $B$, by Lemma \ref{lemmaFactor} we get a birational morphism $\psi:\Sigma\rightarrow X=\pp^2$.
We set $G=\psi_*(\hat B)$, hence $G$ is a curve of degree $10$ which possibly has some [$3,3$]-points and 
no other essential singularities (case ($a$)).
 If $B$ has only [$3,3$]-points as essential singularities, let [$p^\prime\rightarrow p$] be one of them.
Hence, projecting from $p\in \ff_2$ onto 
 the plane we get a birational morphism $\psi:\hat \Sigma \rightarrow \pp^2$ such that $G:=\psi(\hat B)$ is a curve 
 of degree $12$. Moreover, it is easily seen that $G=G^\prime+L_1$, where $L_1$ is a line, and the essential singularities of $G$
are: a $4$-tuple point and a $[5,5]$-point lying on $L_1$, possibly some [$3,3$]-points.
Notice that $L_1$ is the image of the exceptional curve arising from $p$.

Finally, by Lemma \ref{lemmaFactor} $\hat S$ is the canonical resolution of the double
cover of $X$ branched along $G$.      
\end{proof}

\begin{proposition}\label{SII}
Let $S$ be such that $\Hat S$ is a surface of type $\mathcal{S}_{II}$.
Then there exists $n\ge 2$ and  a birational morphism $\psi :\Hat \Sigma \rightarrow \pp^2$ such that
setting $G=\psi_\ast(\Hat B)$ we have:
\begin{itemize}
\item[$a$)]$G=G^\prime +\sum_{i=1}^nL_i,$
where $L_1,\dots ,L_n$ are distinct lines passing through a point $\gamma$ in $\pp^2$ and
 $G^\prime \in|(10+n)L|$ is a reduced curve;
\item[$b$)] the essential singularities of $G$ are a $(2n+2)$-tuple point at $\gamma$, a $[5,5]$-point
$[p^\prime_i\rightarrow p_i]$ such that $p_i\in L_i$, $i=1,\dots, n$,
 possibly some 4-tuple points or $[3,3]$-points.
\end{itemize}
Therefore,
\begin{itemize}
\item[$c$)] $\Hat S$ respectively $S$ are the canonical resolution and the smooth minimal model
of a Du Val double plane of type $\mathcal{D}_n$ with $n\ge 2.$
\end{itemize}
\end{proposition}
\begin{proof}
We have a morphism $\varphi:\Hat \Sigma \rightarrow \ff_e$
such that $\varphi_\ast (\Hat B)=B=B_1+\Gamma_1+\dots + \Gamma_n, \ n\in \{2,\dots ,6\}$, where
$\Gamma_1,\dots ,\Gamma_n$ are pairwise distinct fibres and $B\equiv 8C_0+(4e+8+2(n-1))\Gamma$.

The essential singularities of $B$ are a [5,5]-point $[p_i^\prime\rightarrow p_i]$ such that
$p_i\in \Gamma_i,\ i=1,\dots,n,$ possibly some 4-tuple points or [3,3]-points.

 By Lemma \ref{corollMain} and \ref{lemmaE}, we know that $C_0\not \subset B$ and $0\leq e \leq \frac{6+n}{4}$.
In particular, $e\leq 3$ since $n\leq 6$.

If $e=1$ let $cont_{C_0}:\ff_1\rightarrow \pp^2$ be the birational morphism which contracts $C_0$ to a point
$\gamma$ in $\pp^2$.
We denote by the same letter the image of $p_i$ in $\pp^2$.
Then $L_i:={cont_{C_0}}(\Gamma_i)$ is a line passing through $\gamma$ and
$G^\prime={cont_{C_0}}\sb\ast(B_1)$ is a reduced curve
having an $m$-tuple point at $\gamma$ where
$m=B_1.C_0=n+2$,
a $[4,4]$-point $[p_i^\prime\rightarrow p_i]$ such that $p_i\in L_i$ and $p_i^\prime$ is infinitely near to
$L_i, i=1,\dots,n$, possibly some 4-tuple points or [3,3]-points.

We set $G=G^\prime+L_1+\dots+L_n$ and $\psi=cont_{C_0}\circ \varphi$. Hence a straightforward calculation shows
that $G\in |(10+2n)L|$ and $\psi_\ast(\Hat B)=G$. Whence $(a),(b)$ follow and by Lemma \ref{lemmaFactor} we get $(c)$.

If $e=2$ we can assume that $p_1\not \in C_0$. In fact, if were $p_i\in C_0,$ for $i=1,\dots,n,$ then it would be
$5n \leq C_0.B=-16+16+2(n-1)=2n-2$, a contradiction.

Let us perform the elementary transformation $elm_{p_1}:\ff_2 \dashrightarrow \ff_1$ and consider the curve
\[
B^\prime=B_1^\prime+\Gamma_1^\prime+\dots+\Gamma_{n-1}^\prime+\Gamma_n^\prime
\]
where $\Gamma_i^\prime$ (resp. $B_1^\prime$) is the proper transform of $\Gamma_i,\ i=2,\dots, n$, (resp. $B_1$)
and $\Gamma_1^\prime$ is the (image of the) exceptional curve arising from $p_1$.

Then the proper transform $C_0^\prime$ of $C_0$ is the $(-1)$-section and $C_0^\prime\not \subset B^\prime.$
As we blow up at a singular point of $B$ with odd multiplicity, by Lemma \ref{lemmaFactor}
there exists a birational morphism $\varphi^\prime:\hat \Sigma \rightarrow \ff_1$ such that $\varphi^\prime(\Hat B)=B^\prime$.

Moreover, a straightforward calculation
shows that $B'.C^\prime_0=2n+2$ and $B'$ has the same singularities as $B$.
Therefore we conclude as above.

If $e=3$ we have $n=6$ since $n\geq 4e-6=6$ and so $C_0.B=6$. Hence we see that $p_i\not\in C_0,\ i=1,\dots,6$.
Consider the birational map $elm_{p_1}\circ elm_{p_2}: \ff_3 \dashrightarrow \ff_1$.

Then as above we have a morphism $\varphi^\prime :\hat \Sigma\rightarrow \ff_1$ such that $\varphi^\prime(\Hat B)=B^\prime$,
where $B^\prime=B_1^\prime+\Gamma_1^\prime+\dots +\Gamma^\prime_6$ is composed by the proper transforms of
$B_1,\Gamma_3^\prime,\dots, \Gamma_6^\prime$ and
by the exceptional curves $\Gamma_1^\prime,\Gamma_2^\prime$ arising from $p_1,p_2$.

Also in this case we get $C^\prime_0.B^\prime=2n+2$ and $B^\prime$ has the same singularities as $B$.
Therefore we conclude as above.

If $e=0$ we argue as in the other cases.  
\end{proof}

\section{Du Val double planes}

We are going to complete the proof of Theorem \ref{mainth}. In particular, we will prove implication $(c)\Rightarrow (b)$
and assuming $(c)$ we will show that $S$ is regular unless $p_g(S)=q(S)=1$.

Hence, throughout this section we will assume that $S$ is a minimal surface of general type
which  is the smooth minimal model of a Du Val double plane $X\sp \prime$ and such that 
does not present the standard case.

We will denote by $S^\ast$ the canonical resolution of such a double plane, so we have the following commutative diagram
\[
\xymatrix@!0{
 & & & & &  S\ar@{.>}[dlldd]^\phi \ar@{.>}[lld]  \\
S^\ast \ar[rrr] \ar[rrrrru]^\pi \ar[dd]_{\tilde \rho} &   &    & X^\prime  \ar[dd]_\rho          \\
                \\
W_s \ar[rrr]^\omega   &    &    & W_0&=X
  }
\]
where $X$ is either $\pp^2$ or $\ff_2$ according
to the type of $X^\prime$ (cfr. introduction) and $\rho,\tilde \rho$ are double covers branched
along $G, G_s=\omega^\ast(G)-\sum 2[\frac{m_i}2]\calE_i^\ast$
respectively (cfr. $(1,2)$). 
Furthermore, there is an involution $\sigma^\ast$ on $S^\ast$ induced by $\tilde\rho$ whose fixed locus is the divisor
$R^\ast:=\tilde\rho^{-1}(G_s)$.
We denote by $\Delta\in Pic(X)$ (resp. $\Delta_s\in Pic(W_s)$) a divisor such that $G\in |2\Delta|$ 
(resp. $G_s\in|2\Delta_s|$).

\textbf{Notation.} Let $X^\prime$ be a Du Val double plane of type $\mathcal {D}_n$. We denote by $\delta_1$
the number of $[3,3]$-points of the branch curve $G$, whereas by $\delta_2$ we denote the number of $4$-tuple
points.

Furthermore, if $n>0$ we denote by $[p_i^\prime\rightarrow p_i]$ the $[5,5]$-point of $G$ such that
$p_i\in L_i, i=1,\dots,n$, whereas if $\delta_1>0$ (resp. $\delta_2>0$) we denote by $[q_j^\prime\rightarrow q_j]$
(resp. $r_j$) a $[3,3]$-point (resp. $4$-tuple point) of $G, j=1,\dots,\delta_1$ (resp. $\delta_2)$.

\begin{lemma}\label{can} Let $S^\ast$ be the canonical resolution of a Du Val double plane of type $\mathcal D_n$. Then
\begin{itemize}
\item[$i$)] $p_g(S^\ast)-q(S^\ast)=6-n-\delta_1-\delta_2$
\item[$ii$)] $K_{S^\ast}^2=8-2n-2\delta_1-2\delta_2$
\end{itemize}
Moreover,
\begin{itemize}
\item[$iii$)] $n+\delta_1+\delta_2\le 6$;
\item[$iv$)] if $n\le1$, then $\delta_2\le n$.
\end{itemize}
\end{lemma}
\begin{proof} By \cite{H} we  have
 \begin{align*}
\chi(S^\ast)&=\frac12(K_{\pp^2}+ \Delta ).\Delta+2\chi(\pp^2)-
\frac12\sum\left[\frac{m_i}2\right]\left(\left[\frac{m_i}2\right]-1\right)= \\
&=\frac12 (2+n)\cdot (5+n)+2-\frac12((n+1)n+(2+6)n+2\delta_1+2\delta_2)= \\
&= 7-n-\delta_1-\delta_2
\end{align*}
and 
 \begin{align*}
K_{S^\ast}^2&=2(K_{\pp^2}+ \Delta )^2-2\sum\left(\left[\frac{m_i}2\right]-1\right)^2 =\\
&=2(2+n)^2-2(n^2+(1+4)n+\delta_1+\delta_2)=\\
&=8-2n-2\delta_1-2\delta_2.  
\end{align*}
Since $\chi(S^\ast)=\chi(S)\ge1$, $i)$ implies $n+\delta_1+\delta_2\ge 6$.  

Finally, assume that $n=0$ (resp. $n=1$) and $\delta_2\ge 1$ (resp. $\ge 2$). Let $\mathcal C\subset W_s$
be the strict transform of a general line (resp. conic) passing through
$r_1$ (resp. $r_1,r_2,p_1,p_1^\prime$) and $\tilde{\mathcal C}$ its pull back to $S^\ast$.
Then $|\tilde{\mathcal C}|$ is a pencil of curves of genus 2. A contradiction.    
\end{proof}

Notice that if $n=1$ then $\gamma$ is a $4$-tuple point, hence $\gamma$ may be infinitely near to $p_1^\prime.$

Recall that $\omega$ factors as $\omega_s\circ \dots \circ \omega_1$ where $\omega_{i+1}$ is the blow up of 
$y_i\in W_i$ with exceptional curve $\calE_{i+1}$, $i=0,\dots,s-1$.
\begin{lemma}\label{-2c}
Let $S^\ast$ be the canonical resolution of a Du  Val double plane
of type $\mathcal D_n$ and let $C$ be a reduced and irreducible curve on $W_s$. Then
\begin{itemize}
\item [(1)] $C$ is a $(-2)$-curve contained in $G_s$ such that $\omega(C)=p$ is a point if and only if
there exists $i\in\{0,\dots,s-2\}$ such that $G_i$ has an $[r,r]$-point at $y_i$ with $r\ge 3$ odd,
$(\omega_{i+2}\circ\dots\circ\omega_s)(C)=\calE_{i+1}$ and $(\omega_1\circ \dots\circ\omega_{i-1})(y_i)=p$
(or $y_0=p$ if $i=0$).
\item [(2)]Assume that $n\ge 1$. Then $C$ is a $(-2)$-curve contained in $G_s$ such that $\omega(C)=L$ is a line passing
 through $\gamma$ (resp. $p_i$, $i=1,\dots,n$) if and only if $L\in \{L_1,\dots,L_n \}.$
\end{itemize}
\end{lemma}
\begin{proof}
$(1)$ is straightforward.

$(2)$. If $\omega(C)=L$ is a line then $C$ is the strict transform of $L$, because it is reduced and irreducible.
Notice that since $L_i(G-L_i)= 2n+9$ the only singular points of $G$ lying on, or infinitely near to,
$L_i$ are $\gamma,p_i,p_i^\prime$.

If $L$ is a line passing through $\gamma$ we can assume $y_0=\gamma$ and hence we have that
\[
C=\omega^\ast(L)-\calE_1^\ast-\sum\sb {i\ge 2} c_i\calE_i^\ast
\]
where $c_i=1$ if and only if
$y_{i-1}$ lies on, or is infinitely near to, $L$  and $c_i=0$ otherwise.
Thus, $L^2-C^2=1-1-\sum c_i$ and so $C\sp 2=-2$ if and only if there are exactly two $c_i$'s
which are non zero. Therefore, we get
\[
-2=C.G_s=L.G-(2n+2)-2[\frac{m_j}2]-2[\frac{m_k}2]
\]
where  $j,h\in \{2,\dots,s\}$
are such that $[\frac{m_j}2]+[\frac{m_k}2]=5$.
It is easy to check that the only possibility is $L\in\{L_1\dots,L_n\}$.

If $p_i\in L$ then $L$ is tangent to $L_i$, since $L\subset G$, and so $L=L_i$.
\end{proof}

Now for the proof of the main theorem we consider three cases:
\begin{itemize}
\item[A)] $X^\prime$ is of type $\mathcal D_n$ with $n\ge 2$;
\item[B)] $X^\prime $ is of type $\mathcal D_n$ with $n<2$;
\item[C)] $X^\prime$ is of type $\mathcal D$ or $\mathcal B$. 
\end{itemize}

\subsection*{Case A). $X^\prime$ is of type $\mathcal D_n$, $n\ge 2$.}

As $2n+2\ge 6$ we can assume $y_0=\gamma$. Moreover, since the $[5,5]$-points are the singularities of $G$
 with the greatest multiplicity lower than
$2n+2$ we can assume $y_{i}=p_i$ and $y_{n+i}=p^\prime_i$, $i=1,\dots,n$.
Finally, we assume $y_{2n+i}=r_i$ for $i=1,\dots,\delta_2$ and $y_{2n+\delta_2+j}=q_j$ (resp. 
$y_{2n+\delta_2+\delta_1+j}=q_j^\prime$) for $j=1,\dots,\delta_1$.

\begin{proposition}\label{n>1}
Let $S$ be the smooth minimal model of a Du Val double plane of type $\mathcal{D}_n$ with $n\ge2$.
Then \small{
\begin{itemize}
\item[$a$)]
\begin{align*}
p_g(S)-q(S)&=6-n-\delta_1-\delta_2 \\
K_{S}^2&=8-\delta_1-2\delta_2
\end{align*}
\item[$b$)]there is a rational pencil $|H|$ on $S$ such that:
\begin{itemize}
\item[$i)$] the general member $H\in|H|$ is a smooth hyperelliptic curve of \hbox{genus 3};
\item[$ii)$] $|H|$ has $n$ double curves;
\item[$iii$)] $|H|$ does not have base points.
\end{itemize}
\item[$c$)] The bicanonical map of $S$ factors through $\rho$ and it induces the hyperelliptic involution on the general
$H\in|H|$.
\end{itemize}}
\end{proposition}

\begin{proof}

Let $L$ be a general line in $\pp^2$ passing through $y_0$ and $\tilde L=\omega^\ast(L)-\calE_1^\ast$
its strict transform. Let $\tilde H=\tilde \rho^\ast(\tilde L)$ be the pull back of
$\tilde L$ to $S^\ast$. Therefore, 
$\tilde \rho|_{\tilde H}:\tilde H\rightarrow \tilde L$ is a double cover branched in $\tilde L.G_s=8$ points
and $\tilde H$ is
a smooth hyperelliptic curve of \hbox{genus 3}.
Moreover, $|\tilde H|$ is a rational pencil such that
\[
\tilde H^2=0,\ \ \ \tilde H. K_{S^\ast}=4,\ \ \ \tilde H.R^\ast=8.
\]
For each $i=1,\dots, n$ we set $\tilde L_i=\omega^\ast(L_i)-\calE_1^\ast\in|\tilde L|$ and
we denote respectively by $C_i=\omega^\ast(L_i)-\calE_1^\ast-\calE_{i+1}^\ast-\calE_{i+n+1}^\ast$
the strict transform of $L_i$ and by $C_{n+i}=\calE_{i+1}^\ast-\calE_{i+n+1}^\ast$ the strict transform of
$\calE_{i+1}=\omega_{i+1}^{-1}(y_i)$ on $W_s$.

Hence we have $\tilde L_i=C_i+C_{n+i}+2\calE_{n+1+i}^\ast$ and $C_1,\dots,C_{2n}$ are $(-2)$-curves
 belonging to $G_s$, by Lemma \ref{-2c}.
Therefore, setting $\tilde H_i=\tilde\rho^\ast(\tilde L_i)$ and $E_i=\tilde\rho^{-1}(C_i)$ we get
$\tilde H_i=2E_i+2E_{n+i}+2\tilde\rho^\ast(\calE_{n+1+i}^\ast)$, i.e. $\tilde H_i$ is a double curve.

If $\delta_1>0$, by Lemma \ref{-2c} there are $\delta_1$ more $(-2)$-curves $C_{2n+1},\dots,C_{2n+\delta_1}$ arising from the
[$3,3$]-points which belong to $G_s$.
We set $E_i=\tilde\rho^{-1}(C_i),\ i=2n+1,\dots,2n+\delta_1$.

Notice that the $E_i$'s are $(-1)$-curves on $S^\ast$ and since $G_s$ is smooth they are pairwise disjoint.
 Moreover, it is easily seen that $C_i.\tilde L=E_i.\tilde H=0$.

Since $S$ is minimal of general type,
the birational morphism $\pi:S^\ast \rightarrow S$ factors as $\pi_2\circ \pi_1$ where $\pi_1:S^\ast\rightarrow S^\prime$
contracts (exactly) the $E_i$'s. Hence, by Lemma \ref{can} we get
\[
\begin{array}{l}
p_g(S^\prime)-q(S^\prime)=p_g(S^\ast)-q(S^\ast)=6-n-\delta_1-\delta_2 \\
\\
K_{S^\prime}^2=K_{S^\ast}^2+2n+\delta_1=8-\delta_1-2\delta_2
\end{array}
\]
and so the following table for $(\chi(S^\prime)-1, K_{S^\prime}^2)$:\begin{dimi}
\[
\xymatrix@R-=.7pc@C-=.6pc
{
\chi(S^\prime)-1 &    &   &    &    & & & K_{S^\prime}^2 & & & \\
4   &&&     &    &                    &          &        &                    & 8 \ar[lld] \ar[d]  \\
3   &&&     &    &                    &          & 6\ar[d]\ar[lld]  &    & 7 \ar[lld] \ar[d] &  &  & 8 \ar[lld] \ar[d]  \\
2   &&&     &    &   4\ar[d]  \ar[lld] &          & 5\ar[d]\ar[lld]  &  & 6 \ar[lld] \ar[d] & 6  \ar[d] \ar[lld] |!{[ll];[d]}\hole &  & 7 \ar[lld] \ar[d]& & & 8 \ar[lld] \ar[d]  \\
1  &&& 2  \ar[d] \ar[dll]   &    &     3  \ar[d]   \ar[lld]  &          & 4\ar[d]\ar[lld] &  4  \ar[d] \ar[lld] |!{[ll];[d]}\hole    & 5 \ar[d]\ar[lld] |!{[ll];[d]}\hole & 5 \ar[d]\ar[lld] |!{[ll];[d]}\hole &
 &  6 \ar[lld] \ar[d]  &6 \ar[lld]|!{[ll];[d]}\hole \ar[d] & & 7 \ar[lld] \ar[d]&
8 \ar[lld]|!{[ll];[d]}\hole \ar[d]  \\
0 &0&& 1 & & 2 &  2 & 3&  3 &  4   &  4  &  4   &  5  & 5 & 6& 6 & 7&8&&&&&
  }
\]\end{dimi}
where $K_{S^\prime}\sp 2=8$ if an only if $\delta_1=\delta_2=0$ if and only if $G$
has neither 4-tuple points nor [3,3]-points and the arrow $\swarrow$ (resp. $\downarrow$)
means that one imposes one more 4-tuple point (resp. [3,3]-point) to $G$.

Set $H^\prime=\pi_{1\ast}(\tilde H)$. Then the general member of $|H^\prime|$
is a smooth hyperelliptic curve of genus $3$ because $\tilde H.E_i=0$,
and $H_j^\prime:=\pi_{1\ast}(\tilde H_j)$ is a double curve for each $j=1,\dots,n$.
In particular, we have
\[
{H^\prime}^2=0,\ \ \ K_{S^\prime}.H^\prime=4,\ \ \ R^\prime.H^\prime=8.
\]

Notice that $\sigma^\ast$ induces an involution $\sigma^\prime$ on $S^\prime$ which is a morphism and whose fixed locus
is union of the smooth curve $R^\prime:={\pi_1}_*(R^\ast)$ and the points
$\pi_1(E_1),\dots,\pi_1(E_{2n+\delta_1})$.

Let $H$ denote the image on $S$ of a general $H^\prime\in|H^\prime|$.
Suppose that \hbox{$S^\prime=S$}. Hence $(a)$ and $(b)$ follow. Moreover, as $S$ is minimal and $W_s$
is rational we can apply Proposition \ref{wellknow}.
Therefore, we have
$h^i(W_s,\calO_{W_s}(2K_{W_s}+\Delta_s))=0$, $i>0,$ and for $(c)$ it suffices to
show that $h^0(W_s,\calO_{W_s}(2K_{W_s}+\Delta_s))=0$.

On the other hand by the Riemann-Roch formula we get
\[
\begin{split}
h^0(W_s,&\calO_{W_s}(2K_{W_s}+\Delta_s))=\chi(2K_{W_s}+\Delta_s)= \\
&=\frac12(2K_{\pp^2}+\Delta).(K_{\pp^2}+\Delta)-\frac18\sum(m_i-4)(m_i-2)+1 \\
&=\frac12(n^2+n-2)-\frac18(4n^2+4n)+1=0
\end{split}
\]
whence the bicanonical map of $S$ factors through $\rho$.

So it remains to prove that $S^\prime=S$. Suppose to the contrary that $S^\prime\neq S$, then
$\pi_2:S^\prime\rightarrow S$ is not the identity and there is a $(-1)$-curve $E\subset S^\prime$ contracted by $\pi_2$
 to a point.

First of all we claim that $E.H^\prime=0$. In fact, $E.H^\prime\ge 0$ since $|H^\prime|$ is a pencil. 

If $E.H^\prime>0$ then $E.H^\prime\ge 2$ since $|H^\prime|$ has $n\ge2$ double curves. Then we get
$H^2\ge 4,H.K_S\le 2$ and so the Hodge Index
theorem implies that $K_S^2=1$ and $H$ is numerically equivalent to $2K_S$.
Observe that in this case $X^\prime$ is of type $\mathcal D_2$ with $\delta_1=0, \delta_2=4$. In particular,
the involution $\sigma^\prime$ acting on $S^\prime$ has 4 isolated fixed points.

As $K_S^2=1$ then $\pi_2$ contracts exactly $E$, and $H^2=4$ implies $E.H^\prime=2$.
Now we have to consider two cases: either $E$ belongs to $R^\prime$ or not.

If $E$ belongs to $R^\prime$ then $q:=\pi_2(E)$ is an isolated 
fixed point of the induced involution $\sigma$ on $S$ since $R^\prime$ is smooth. Moreover we get 
$K_S.R=\frac12H.R=3$
where we denote by $R:={\pi_2}_\ast(R^\prime)$ the divisorial part of $Fix(\sigma)$. A contradiction, indeed $\sigma$
has $5$ isolated fixed points and so by Proposition \ref{wellknow} we get $K_S.R=1$.

If $E$ does not belong to $R^\prime$ then $2K_S.R=H.R\ge 8$.
On the other hand, in this case $\sigma$ has $4$ isolated fixed points and so by Proposition \ref{wellknow}
it follows $4=K_S.R+4$. A contradiction.

Therefore, $H^\prime.E=0$.

Let $\st E\in S^*$ be the strict transform of $E$, we set $\calE=\tilde \rho (\st E)$.
Therefore, we have $\calE.\tilde L=\st E.\tilde H=E.H^\prime=0$ and then $\calE$ is a component
of a curve $\tilde L_E\in|\tilde L|$. In particular, $\calE$ is a smooth rational curve and
$\tilde \rho|_{\st E}:\st E\rightarrow \calE$ is either
a double cover or an isomorphism.
We consider the two cases separately.

\textbf{If $\tilde \rho|_{\st E}$ is an isomorphism.} We have $\tilde\rho^\ast(\calE)=\st E+\tilde E$ where
$\st E \cong \tilde E$ (possibly $\st E =\tilde E$).

If $\st E=\tilde E$ then $\calE\subset G_s$ and thus $\st E\cap E_i=\emptyset, i=1,\dots,2n+\delta_1$
since $G_s$ is smooth.
Hence $-1=E^2=\frac12\calE^2$. By Lemma \ref{-2c} we get a contradiction.

Therefore, $\st E\neq \tilde E$. In this case we have $G_s|_{\calE}=2z$ where $z\in Pic(\calE)$ and 
$\st E^2=\tilde E^2=\calE^2-deg(z)$.
Then either $G_s\cap\calE=\emptyset $ or $G_s$ and $\calE$ are tangent at each intersection point.
In particular $E\cap E_j=\calE\cap C_j=\emptyset, j=1,\dots,2n+\delta_1$ since $\calE\cap C_j\neq \emptyset$ implies
that both $\calE$ and $C_j$ belong to $\tilde L_E$.
Hence $\st E,\tilde E$ are $(-1)$-curves and either $\calE^2=0$ or
$deg(z)=0$, since $\calE^2\le 0$ because $\calE\subset \tilde L_E\in|\tilde L|$.
If $\calE^2=0$ then $\calE.G_s=2$ and $\tilde L_E=a\calE$ for some $a\ge1$.
Hence $a=4$ since $\tilde L.G_s=8$. A contradiction, since $|\tilde L|$ does not have multiple curves.

So $deg(z)=0$ and $\calE^2=-1$. By the definition of canonical resolution $L_E:=\omega(\calE)$ can not be a point,
therefore $L_E$ is a line passing through $\gamma$ and $L_E\neq L_j, j=1,\dots n$.
Since $L_E^2-\calE^2=2$ there is
exactly one point $y_i\neq \gamma$ lying on $L_E$. Analogously to Lemma \ref{-2c} we get that $y_i$
is an $8$-tuple point \hbox{of $G$}. A contradiction.

\textbf{If $\tilde\rho|_{\st E}$ is a double cover.} Then $\st E^2=2\calE^2<0$ is even and hence $\st E^2\le-2$.
Therefore, $\st E.(E_1+\dots+E_{2n+\delta_1})=-1-\st E^2$ is an odd (non zero) number and it is equal to the number of
$E_i$'s which meet $\st E$ since $E$ is smooth.
By the Hurwitz formula we have 
$\st E.(E_1+\dots+E_{2n+\delta_1})\le \st E.R^\ast=\calE.G_s=2$ which yields $\st E^2=-2$ and $\calE^2=-1$.
The usual calculation shows that $\omega(\calE)$ can not be neither a point nor a line through $\gamma$.
Whence $ \Hat S=S$ and the claim follows.
\end{proof}

\begin{proposition}\label{q}
Let $S$ be the smooth minimal model of a Du Val double plane of type $\mathcal{D}_n$ with $n\ge2$.

Then $q(S)=0$ unless $q(S)=p_g(S)=1$.
More precisely, let $\mathcal P$ be the set of $n+\delta_1+\delta_2$ points
$ \{p_1,\dots,p_n,q_1,\dots,q_{\delta_1},r_1\dots,r_{\delta_2}\}$. Then $p_g(S)=q(S)=1$ if and only if
$n+\delta_1+\delta_2=6$ and
\begin{itemize}
\item[-] either no point of $\mathcal P$ is infinitely near to $\gamma$ and the points of $\mathcal P$ lie on a conic;
\item[-] or exactly a point $p\in \mathcal P$ is infinitely near to $\gamma$ and there is a conic passing 
through the set of points $\{\gamma\} \cup \mathcal P\setminus \{p\}$.
\end{itemize}
\end{proposition}
\begin{proof}
Recall that
\begin{align*}
p_g(S)=p_g(S^\ast)&= h^0(W_s,\calO_{W_s}(K_{W_s}+\Delta_s))+h^0(W_s,\calO_{W_s}(K_{W_s}))=\\
&=h^0(W_s,\calO_{W_s}(K_{W_s}+\Delta_s))
\end{align*}
where $\Delta_s\in Pic(W_s)$ is such that $G_s\equiv 2\Delta_s$.

Suppose that $q(S^\ast)\neq 0$, hence $p_g(S^\ast)\neq 0$ and there exists a curve $C\in |K_{W_s}+\Delta_s|$. We have 
\begin{align*}
C & \equiv  \omega^\ast(K_{\pp^2}+\Delta)-\sum ( \left[ \frac{m_i}2 \right] -1)\calE_i^\ast\equiv \\
&\equiv \omega^\ast((2+n)l)-n\calE_1^\ast-\sum_{i=2}^{n+1}(\calE_i^\ast+2\calE_{n+i}^\ast)-\sum_{i=2n+2}^{2n+\delta_2+1}\calE_i^\ast-
\sum_{i=2n+2+\delta_2+\delta_1}^{2n+\delta_2+2\delta_2+1}\calE_i^\ast
\end{align*}
where $l$ is a line in $\pp^2$.

On the other hand, using the notations introduced before, we get the following equalities
\begin{align*}
C&.C_i=C.(\omega^\ast(l)-\calE_1^\ast-\calE_{i+1}^\ast-\calE_{i+n+1}^\ast)=-1,\ i=1,\dots,n; \\
C&.C_{i+n}=C.(\calE_{i+1}^\ast-\calE_{i+n+1}^\ast)=-1,\ i=1,\dots,n; \\
C&.C_{i+2n}=C.(\calE_{i+2n+\delta_2}^\ast-\calE_{1+2n+\delta_2+\delta_1}^\ast)=-1,\ i=1,\dots,\delta_1
\end{align*}
which imply that the $C_i$'s are fixed components of $|C|$. Therefore, we can write
\[
|C|=|\omega^\ast(2l)-\sum_{i=2}^{n+1}\calE_{n+i}^\ast-
\sum_{i=2n+2}^{2n+\delta_2+1}\calE_i^\ast-\sum_{i=2n+2+\delta_2}^{2n+\delta_2+\delta_2+1}\calE_i^\ast|
+\sum_{i=1}^{2n+\delta_1} C_i
\]
and so
\[
p_g(S^\ast)=h^0(W_s, \calO_{W_s}(\mathcal C))
\]
where 
\[
\mathcal C\in |\omega^*(2l)-\sum_{i=1}^{n} \omega^*(p_i)-\sum_{i=1}^{\delta_1}\omega^*(q_i)-
\sum_{i=1}^{\delta_2}\omega^*(r_i)|
\]

Now there are two cases to be considered: either at least a point of the set $\mathcal P=\{p_1,\dots,p_n,
q_1,\dots,q_{\delta_1},r_1,\dots,r_{\delta_2}\}$ is infinitely near to $\gamma=y_0$, or not.

First we consider the second case. Then $p_g(S^*)$ is equal to the dimension of the vector space
$V_m\subset H^0(\pp^2,\calO_{\pp^2}(2l))$ consisting of 
of those conics in $\pp^2$ which passe through the $m= n+\delta_1+\delta_2\le 6$ points of $\mathcal P$.

It is well known that the dimension of  $V_m$ is greater than or equal to $6-m$ and 
by Proposition \ref{n>1} $q(S^*)=0$ if and only if the equality holds. In particular, we can assume $m>3$.

If $4\le m\le 5$ then $p_g(S^*)>6-m$ if and only if there exists a line passing through at least $4$ points
of $\mathcal P$.
Whereas if
$m=6$ then $p_g(S^*)>0$ if and only if all the points lie on a conic and $p_g(S^*)>1$ if and only if at least 
$5$ points are contained in a line.

Assume that $n+\delta_1+\delta_2=4$ and suppose that there exists a line $L^\prime$ passing through the four points.
If $L^\prime \not \subset G$ we get $10+2n = L^\prime.G\ge 5n+3(4-n)=2n+12$, a contradiction.
 On the other hand if $L^\prime \subset G$ then $L^\prime$ is tangent to $G$ at each [$5,5$]-point
  (resp. [$3,3$]-point) and hence $10+2n-1= (G-L^\prime).L^\prime=8n+3(4-n)=5n+12$, a contradiction.

Now assume that $n+\delta_1+\delta_2=5$ and suppose that there exists a line $L^\prime$ passing through four points
of $\mathcal P$.
In particular,
there is at most one $[5,5]$-point which does not lie on $L^\prime$. Hence, either
$10+2n=G.L^\prime\ge 5(n-1)+1+3(4-(n-1))=2n+11$ or
$9+2n=(G-L^\prime).L^\prime\ge 8(n-1)+3(4-(n-1))=5n+7$ depending on $L^\prime \not \subset G$ or $L^\prime \subset G$.
A contradiction.

Finally, assume that $n+\delta_1+\delta_2=6$ and suppose that a line $L^\prime$ passes through $5$ of the points
points. Then there is at most one $[5,5]$-point which does not lie on $L^\prime$.
Hence, either
$10+2n=G.L^\prime\ge 5(n-1)+1+3(5-(n-1))=2n+14$ or
$9+2n=(G-L^\prime).L^\prime\ge 8(n-1)+3(5-(n-1))=5n+10$ depending on $L^\prime \not \subset G$ or $L^\prime \subset G$.
A contradiction.

Therefore, $p_g(S^\ast)>6-n-\delta_1-\delta_2$ implies that $n+\delta_1+\delta_2=6$, the points of $\mathcal P$
lie on a conic
and no five of them are collinear. Hence $\chi(S)=1$ and $p_g(S^*)=q(S^*)=1$.  

Next we discuss the other case. Let us denote by $\st\calE_1\subset W_s$ the strict transform of $\calE_1$.
 First suppose that exactly a point, say $p\in \mathcal P$, is infinitely near to $y_0=\gamma$.
 Hence 
\[
\mathcal C.\st\calE_1\le\mathcal C.(\calE_1^\ast-\sum\sb{i\ge 2} c_i\calE_i^\ast\le\mathcal C.(\calE_1^*-\calE^*)=-1
\]
where $\calE^*=\omega^*(p)$ and $c_i$ is equal to $1$ or $0$ depending on
$y_{i-1}$ is infinitely near to $y_0$ or not.
Therefore, $\st\calE_1$ is a fixed component of $|\mathcal C|$ and $p_g(S^\ast)$ is easily
seen to be equal to the dimension of the vector space consisting of those conics passing
through the set of points $\{\gamma\}\cup\mathcal P\setminus \{p\}$. As before we get the claim. 

Now suppose that at least two points of $\mathcal P$ are infinitely near to $y_0$ and let $p,q$ be two of them.
Let $\st L_i\subset W_1$ be the strict transform of $L_i$ under $\omega_1$, $i=1,\dots,n,$ and denote by
$mult_p(G_1-\sum \st L_i), mult_q(G_1-\sum \st L_i)$ the multiplicity of $G_1-\sum L_i$ at $p$ and $q$, respectively.
Then from the inequalities
\[
6\le mult_p(G_1-\sum L_i)+mult_q(G_1-\sum L_i)\le (G_1-\sum \st L_i).\calE_1=2+n
\]
it follows that $n\ge 4$ and it is easy to check that one has
$n+\delta_1+\delta_2=6$, where $n\in \{4,5,6\}$, and that there are 
exactly two points infinitely near to $\gamma$ which have to be  respectively $p_1,p_2;\ p_1,q_1;\ q_1,q_2$.

We consider the case $n=6$, the others are completely analogous.
Then we have $\mathcal C.\st\calE_1=\mathcal C.(\calE_1^*-\calE_2^*-\calE_3^*)=-2$
and $|\mathcal C|=|\mathcal C^\prime|+\st\calE_1$ where
$\mathcal C^\prime$ is strict transform of a conic through $\gamma, p_3,\dots,p_6$.
Hence $1\le p_g(S^*)=q(S^*)\le 2$ and $p_s(S^*)=2$ if and only if $p_3,\dots,p_6$ lie on a line.

Suppose that $p_g(S^*)=2$. Let $L^\prime$ be the line passing 
through $p_3,\dots,p_6$ and consider the linear system 
$|\omega^*(4l)- 2\calE_1^*-\sum_{i=1}^6(\calE_{i+1}^*+\calE_{i+n+1}^*)|$.
Let $\st F_1,\st F_2$ be the strict 
transforms on $W_s$ of $F_1:=L_1+L_2+2L^\prime$ and \hbox{$F_2:=L_3+\dots+L_6$}, respectively.
Then, $\st F_1,\st F_2$ do not have common components and
 $\st F_j\in|\omega^*(4l)- 2\calE_1^*-\sum_{i=1}^6(\calE_{i+1}^*+\calE_{i+n+1}^*)|$, $j=1,2$.
 Arguing as in \hbox{Lemma \ref{lemma}} we get that the general 
element $F\in |\omega^*(4l)- 2\calE_1^*-\sum_{i=1}^6(\calE_{i+1}^*+\calE_{i+n+1}^*)|$ is a smooth curve of genus 2 such 
that $F.G_s=0$. A contradiction, because we are assuming that $S$ does not present the standard case.             
\end{proof}

\begin{remark} If $n=1$ and the $4$-tuple point lying on $L_1$ is not infinitely near to the [$5,5$]-point, the 
above theorems holds also for $S^\ast$ of type $\mathcal D_1$.
\end{remark}
\subsection*{Case B).  $X^\prime$ is of type $\mathcal D_n$, $n<2$.}

\begin{proposition}\label{n<2}
Let $S$ be the minimal model a Du Val double plane of type $\mathcal{D}_n$ with $n\le 1$.
Then
\begin{itemize}
\item[$a$)]
\begin{align*}
K_{S}^2&=8-n-\delta_1-2\delta_2 \\
p_g(S)-q(S)&=6-n-\delta_1-\delta_2 \\
q(S)&=0\ unless\ p_g(S)=q(S)=1\ and\ K_S^2=3;
\end{align*}
in particular, $p_g(S)=q(S)=1$ if and only if $n=1,\delta_1=5,\delta_2=0$ and the points $p_1,q_1,\dots,q_5$ lie on
a conic;
\item[$b$)] the bicanonical map of $S$ factors through $\rho$; 
\item[$c$)] either $p_g(S)=6, K_S^2=8$ and $K_S$ is ample or there is a rational pencil
$|H|$ on $S$ such that:
\begin{itemize}
\item[$i$)]  the general member $H\in|H|$ is a smooth hyperelliptic curve of genus 3;
\item[$ii$)] the bicanonical map of $S$ induces the hyperelliptic involution on the general
$H\in|H|$.
\item[$iii$)] either $p_g(S)=6, K_S^2=8$ ($K_S$ is not ample) and $|H|$ does not have base points or $|H|$
 has one base point .
\end{itemize}
\end{itemize}
\end{proposition}

\begin{proof} If $n= \delta_1=0$ (and then $\delta_2=0$) and $G$ is smooth, then $S=S^\ast$
and it is easily seen that $(a),(b)$ hold.
 In particular, $p_g=6, K_S^2=8$ and $K_S$ is ample.

 Hence, we can assume that \begin{itemize}
\item[-] $y_0=p$ if $n=\delta_1=0$;
\item[-] $y_0=q_1$ if $n=0, \delta_1\ge 1$;
\item[-] $y_0=p_1$ if $n=1$.
 \end{itemize}
where in the first case $p\in \pp^2$ is a (non essential) singular point of $G$.

Now we proceed as in the proof of Proposition \ref{n>1}.
Let $\tilde H$ be the pull back to $S^\ast$ of a general line passing through $y_0$. Hence, $|\tilde H|$ is a 
pencil of (smooth) hyperelliptic curves of genus $3$. In particular, 
\[
\tilde H^2=0,\ \ \tilde H.K_{S^\ast}=4,\ \ \ R^\ast.\tilde H=8.
\]

Let $\pi_1:S^\ast\rightarrow S^\prime$ be the birational morphism which contracts the $(-1)$-curves 
$E_1,\dots, E_{\delta_1+2n}$ arising from the [$r,r$]-points, $r=3,5$, and from $L_1$ 
(resp. $\pi_1=id$ if $\delta_1+n=0$).  We set $H^\prime=\pi_1(\tilde H)$.

Note that if $\delta_1+n>0$ we can assume
$E_1=\tilde\rho^{-1}(\st\calE_1)$ where $\st \calE_1\subset G_s$ is the strict transform of $\calE_1=\omega_1^{-1}(y_0)$.
Hence, 
$E_i.\tilde H=0$ for each $i>1$ and
\[
{H^\prime}^2=E_1.\tilde H=
\left\{
\begin{array}{l}
1\ if\ n+\delta_1>0 \\
0\ if\ n+\delta_1=0
\end{array}
\right. 
\]

By Lemma \ref{can} 
$p_g(S^\prime)-q(S^\prime)=6-\delta_1-\delta_2-n$ and we get the following table for 
$(\chi(S^\prime)-1,K_{S^\prime}^2)$:
\[
\xymatrix@!0{
\chi(S)-1 &  & K_{S}^2 & & & \\
6   &  & 8 \ar[d]  &  \\
5   &  & 7 \ar[d]  &  8\ar[lld]|!{[ll];[d]}\hole  \ar[d]  \\
4   &  6 \ar[d]& 6\ar[d] &  7 \ar[lld]|!{[ll];[d]}\hole \ar[d]  \\
3   &  5\ar[d]& 5\ar[d] &  6\ar[lld]|!{[ll];[d]}\hole  \ar[d]  \\
2   &  4\ar[d]& 4\ar[d] &  5 \ar[lld]|!{[ll];[d]}\hole \ar[d]  \\
1  &   3\ar[d] & 3\ar[d]&  4 \ar[lld]|!{[ll];[d]}\hole \ar[d]  \\
0  &   2       & 2 &  3 &
  }
\]

As before it suffices to show that $S^\prime=S$. In particular, recall that as
$S$ does not present the standard case, if $q(S)>0$ then $K_S^2>2\chi(S)$ (cfr. \cite{X1}). 

Let us suppose that $E\subset S^\prime$ is a $(-1)$-curve and define $H=\pi_2(H^\prime)$, where $\pi_2\circ\pi_1=\pi$.
First we prove that $E.H^\prime=0$.

Suppose that $E.H^\prime>0$.

If $n+\delta_1>0$ the Hodge Index Theorem gives $K_S^2\le2$. A contradiction, since
 $K_S^2\ge K_{S^\prime}^2+1\ge3$.

If $n+\delta_1=0$, then $S^\prime=S^\ast$ and $K_{S^\prime}^2=8$. Since $H^2\ge 1$ and $H.K_S\le 3$ it follows
from the Hodge Index Theorem that $K_S^2=9$ and $K_S$ is numerically equivalent to $3H$.
In particular $\pi=\pi_1$ contracts exactly $E$ and $E.\tilde H=1$.

If $E\not \subset R^\ast$, then the involution
$\sigma$ induced on $S$ has no isolated fixed points and by Proposition \ref{wellknow} we get $0=R.K_S-20$.
A contradiction, since we have that
$R.K_S=3R.H\ge 3R^\ast .\tilde H= 24$.

Whence $E\subset R^\ast $ and $\calE:=\tilde \rho (E)$ is a $(-2)$-curve belonging to $G_s$.
By Lemma \ref{-2c} $\calE^\prime:=(\omega_2\circ\dots\circ\omega_s)(\calE)$ is a curve on $W_1$.
Denote by $\st L:=\omega_1^\ast(L)-\calE_1$ the strict transform on $W_1$ of
a general line passing through $y_0$. Since $S^\ast \rightarrow W_1$ 
is a (finite) double cover in a neighborhood of $\st L$, we have $\st L.\calE^\prime=\tilde H.E=1$
and so $\calE^\prime$ is smooth.
Therefore, $\omega_1(\calE^\prime)\subset \pp^2$ is a reduced and
irreducible curve of degree $d\ge 1$ with multiplicity
$d-1=\calE^\prime.\calE_1$ at $y_0$ and smooth elsewhere.

It follows that ${\calE^\prime}^2=2d-1\ge 1 $ and we can assume $y_1\in \calE^\prime.$
On the other hand $G_1$ does not have essential singularities and hence
$(\omega_2\circ\dots\circ\omega_t\circ\tilde \rho)^\ast(y_1)$ is a $(-2)$-cycle on $S^\ast$ (cfr. \cite{BPV}).
In particular, there is
a $(-2)$-curve $E^\prime\subset S^\ast$ such that $E.E^\prime=1$ and hence $\pi(E^\prime)\subset S$ is a
$(-1)$-curve. A contradiction.

Therefore, $E.H^\prime=0$. Now arguing as in the proof of Proposition \ref{n>1} we get a contradiction and then
$S^\prime=S$.   
\end{proof}

\begin{proposition}
Let $S^\ast$ be the canonical resolution of a Du Val double plane of type $\mathcal D_0$ such that $\delta_1\ge 2$.
Let $l_i$ be the line tangent to $G$ at $q_i, i=1,\dots,\delta_1$. Then 
$l_j\neq l_k$ for some $j\neq k$ if and only if
$S^\ast$ is the canonical resolution of a Du Val double plane of type $\mathcal D_1$
with a $4$-tuple point and $\delta_1-2$ $[3,3]$-points.
\end{proposition}

\begin{proof}
Assume that $S^\ast$ is the canonical resolution of a Du Val double plane of type $\mathcal D_0$
and suppose that [$q_1^\prime\rightarrow q_1$], [$q_2^\prime\rightarrow q_2$] are 
such that $l_1\neq l_2$. Then we can perform the quadratic transformation of the plane
$\lambda_{q_1,q^\prime_1, q_2}: \pp^2\dashrightarrow \pp^2$
centered at $q_1,q^\prime_1, q_2$.
Let $G^\prime$ be the proper transform of $G$ under $\lambda_{q_1,q^\prime_1, q_2}$ 
and let $L_1$ be the image of the exceptional curve arising from $q_1$.
Then it is easily seen that $G^\prime$ is a reduced curve of degree $11$ with a triple point and a
[$4,4$]-point lying on the line $L_1$, a $4$-tuple point at the image of $q_1^\prime$
and $\delta_1-2$ [$3,3$]-points at the image of $q_3,\dots,q_{\delta_1}$.

Therefore, we have the following commutative diagram
\[
\xymatrix@!0{
 & & &\ \ S^\ast \ar[ddll]_\omega  \ar[ddrr]^{\omega^\prime}      \\
                \\
G \subset & \pp^2  \ar@{.>}[rrrr]^{\lambda_{q_1,q^\prime_1, q_2}} &    &  &  & \pp^2&\ \ \supset G^\prime+L_1 
  }
\]
where $\omega^\prime$ is a morphism because we blow up singular points of $G$.
Now arguing as in Lemma \ref{lemmaFactor} one sees that $S^\ast$ is the canonical resolution of the double cover
of $\pp^2$ branched along $G^\prime+L_1$.

For the converse, perform the quadratic transformation of $\pp^2$ centered at $p_1,p^\prime_1,r_1$.
\end{proof}

\subsection*{Case C). $X^\prime$ is of type $\mathcal B$ or $\mathcal D$.}

 This is the easiest case. In fact, arguing as above one gets the following:
\begin{proposition}\label{easy}
Let $S$ be the minimal model a Du Val double plane of type $\mathcal B$ or $\mathcal D$.
Then $q(S)=0$ and the bicanonical map of $S$ factors through $\rho$. Moreover, 
\begin{itemize}
\item[$a$)] if $S$ is of type $\mathcal D$, then $p_g(S)=3, K_S^2=2$ and $K_S$ is ample;\\
\item[$b$)] if $S$ is of type $\mathcal B$, then $p_g(S)=6, K_S^2=9$ and there is a rational pencil
$|H|$ on $S$ such that:
\begin{itemize}
\item[$i$)] the general member $H\in|H|$ is a smooth hyperelliptic curve of genus 3;
\item[$ii$)] the bicanonical map of $S$ induces the hyperelliptic involution on the general
$H\in|H|$.
\item[$iii$)] $|H|$ has one base point .
\end{itemize}
\end{itemize}
\end{proposition}

\section{Conclusion and Remarks}

We collect some corollaries of the main theorem. Throughout the end $S$ will be a minimal surface of general type 
not presenting the $standard$ $case$.

\begin{corollary}
Let $S$ be the smooth minimal model of a Du Val double plane of type $\mathcal D_n$. If $n\ge 2 $ then 
$(\mathbb Z_2)\sp {n-1}\subseteq Tors(S)$.
\end{corollary}
\begin{proof}
Because the rational pencil $|H|$ has $n$ pairwise distinct double curves 
the claim is clear. 
\end{proof}

As we remarked in the introduction, if $S$ is the smooth minimal model of a double plane with $p_g(S)\ge 2$ then 
the bicanonical map of $S$ has degree 2, because $S$ is regular.
In the following corollary we show that if $p_g(S)\le 1$ then $\varphi_{2K}$ may have degree greater then 2.
\begin{corollary}
Let $S$ be the smooth minimal model of a Du Val double plane with
$p_g(S)=1,q(S)=0$ and $K_S^2=2$. Then $\varphi_{2K}$ has degree 4 and $S_2$ is a quadric cone in $\pp^3$.
\end{corollary}
\begin{proof}
By Proposition \ref{n>1}, Proposition \ref{q} and Proposition \ref{n<2}, $S^\ast$ is of type $\mathcal D_2$
and the branch curve $G$ has $\delta_1=0$ [$3,3$]-points and $\delta_2=3$ $4$-tuple points. 

Although the claim follows from a result of Catanese, Debarre (cfr. \cite{CD}, Proposition 1.5)
and the previous Corollary, it can also be proved with the same argument we used before.

As the bicanonical map of $S$ factors through the involution induced by the double cover, we have the commutative
diagram
\[
 \xymatrix@!0{
&&S\sp \ast \ar [dd]_{\tilde \rho} \ar[rr]^\pi & & S \ar [dd]^{\varphi\sb{2K}}  \\
 \\
\pp^2&&\hat W\sb s \ar[ll]_\omega \ar[rr]^{\varphi_F}& & S_2
   }
\]
where $\varphi\sb F$ is the morphism defined by the linear system
\begin{equation*}
\begin{split}
|F|=&|2K\sb{W\sb s}+G\sb s-\sum\sb 1\sp 4 C_i|=\\
=&|\omega\sp\ast(6l)-\omega\sp \ast(2\gamma)-\sum_{i=1}^2\omega\sp\ast(2p_i+2p^\prime_i)-\sum_{i=1}^3\omega\sp\ast(2r_i)|
\end{split}
\end{equation*}
(here $l$ is a line in $\pp^2$).
Let $\mathcal C_3$ denote the strict transform under $\omega$ of a general cubic in $\pp^2$ passing through the set
$\mathcal P=\{\gamma,p_1,p_1^\prime,p_2,p_2^\prime,r_1,r_2,,r_3\}$.
Then $\mathcal C_3$ is a smooth curve of genus 1
and the linear system $|\mathcal C_3|$ has one base point $p\in\hat W$.

Now one sees that $|F|$ cuts a $g^1_2$ on the general curve in $|\mathcal C_3|$ and so $\varphi\sb F$ has
degree greater than 2. It follows that $\varphi\sb{2K}$ has degree $d\ge 4$ and $S_2$ is a surface of degree
$\frac8d$ in $\pp^3$. Whence, $d=4$ and $\varphi\sb {2K}(S)=\varphi\sb F(W_s)$ is a quadric cone with
vertex $\varphi\sb F(p)$ and ruling $\varphi \sb F(\mathcal C_3)$.
\end{proof}

\begin{proposition}
Let $S$ be a smooth minimal surface of general type with $p_g(S)=0$. If the bicanonical map has degree 2 and $S$
does not present the standard case, then
\begin{itemize}
\item[1)] either $K_S^2=3$ and $S_2\subset \pp^2$ is an Enriques surface,
\item[2)] or $S$ is the smooth minimal model of a Du Val double plane of \hbox{type $\mathcal D_n$}
with $K_S^2$ and $n$ as in the following table:
\[
\begin{tabular}{|c|c|c|c|c|c| c |c|}
\cline{1-8}
$K_S^2$ & 2       & 3     & 4     & 5   & 6   & 7 & 8 \\
\cline{1-8}
$n$     & 0,1,2,3 & 1,2,3 & 2,3,4 & 3,4 & 4,5 & 5 & 6  \\
\cline{1-8}
\end{tabular}
\]
\end{itemize}

Moreover, in case $2)$ there is a rational pencil $|H|$ whose general member is a smooth hyperelliptic curve of genus 3
such that
\begin{itemize}
\item[-] the bicanonical map of $S$ induces the hyperelliptic involution on the general curve $H\in|H|$;
\item[-] if $n\le 1$ then $|H|$ has one base point; 
\item[-] if $n\ge 2$ then $|H|$ is base points free and has $n$ double fibres.
\end{itemize}
\end{proposition}
\begin{proof}
By \cite{X2} if $S_2$ is not rational then $K_S^2=3,4$ and $S_2$ is an Enriques surface.
In \cite{MP1} M.Mendes Lopes and R.Pardini
show that the case $K_S^2=4$ does not occur.

Now the claim follows by Theorem \ref{mainth}, Proposition \ref{n>1}, Proposition \ref{n<2} and Proposition \ref{easy}.   
\end{proof}

\begin{remark} We remark that the above result was partially proved by R.Pardini and M.Mendes Lopes.
In fact, they classify 
surfaces of general type with $6\le p_g \le 8$ and bicanonical map of degree two in \cite{MP2},\cite{MP3},\cite{P}
where they also construct examples of such surfaces.
\end{remark}

As we remarked in the introduction, we get an analogous result for regular surfaces with $p_g=1$.
\begin{proposition}
Let $S$ be a smooth minimal surface of general type with $q(S)=0$ and $p_g(S)=1$.
If the bicanonical map has degree 2, $S$
does not present the standard case and the bicanonical image $S_2$ is not a K3 surface, then
$S$ is the smooth minimal model of a Du Val double plane of \hbox{type $\mathcal D_n$}
with $K_S^2$ and $n$ as in the following table:
\[
\begin{tabular}{|c|c|c|c|c|c| c |c|}
\cline{1-8}
$K_S^2$ & 2       & 3     & 4     & 5   & 6   & 7 & 8 \\
\cline{1-8}
$n$     &     2   & 0,1,2 & 1,2,3 & 2,3 & 3,4 & 4 & 5  \\
\cline{1-8}
\end{tabular}
\]

Moreover, there is a rational pencil $|H|$ whose general member is a smooth hyperelliptic curve of genus 3
such that
\begin{itemize}
\item[-] the bicanonical map induces the hyperelliptic involution on the general $H\in|H|$;
\item[-] if $n\le 1$ then $|H|$ has one base point; 
\item[-] if $n\ge 2$ then $|H|$ is base points free and has $n$ double fibres.
\end{itemize}
\end{proposition}

Finally, we get a partial result concerning the case $p_g(S)=q(S)=1$.
\begin{proposition}
Let $S$ be a smooth minimal surface of general type with $p_g(S)=q(S)=1$ and $7\le K_S^2\le 8$. Assume
that the bicanonical map of $S$ has degree 2 and that $S$ does not present the standard case. Then 
\begin{itemize}
\item [-] If $K_S^2=7$, then $S$ is the smooth minimal model of a double plane branched along a reduced curve 
$G=G^\prime+L_1+\dots+L_5$, where $G^\prime$ has degree $15$ and $L_1,\dots,L_5$ are lines meeting at a point $\gamma.$ 
The essential singularities of $G$ are a $12$-tuple point at $\gamma$, a $[5,5]$-point $[p_i^\prime\rightarrow p_i]$
on $L_i, i=1,\dots,5$, a $[3,3]$-point $[q_1^\prime\rightarrow q_1]$. The points $p_1,\dots,p_5,q_1$ lie on a conic.
\item[-] If $K_S^2=8$, then $S$ is the smooth minimal model of a double plane branched along a reduced curve 
$G=G^\prime+L_1+\dots+L_6$, where $G^\prime$ has degree $16$ and $L_1,\dots,L_6$ are lines meeting at a point $\gamma.$ 
The essential singularities of $G$ are a $14$-tuple point at $\gamma$, a $[5,5]$-point $[p_i^\prime\rightarrow p_i]$
on $L_i, i=1,\dots,6$. The points $p_1,\dots,p_6$ lie on a conic.
\end{itemize}
\end{proposition}
\begin{proof}
By \cite{X2}, Theorem 3, the bicanonical image $S_2$ is a rational surface. Therefore, we can apply Theorem \ref{mainth}
and then the results of Section 4.    
\end{proof}
\textbf{Remark.} Surfaces with $p_g=q=1,K\sp 2=8$ and bicanonical map of \hbox{degree 2} are studied in detail
and classified by F.Polizzi in his PhD thesis (cfr.\cite{Po}).
In particular, he constructs such surfaces as the quotient of the product of two curves by a finite group. 


\begin{thebibliography}{Hor}
\small {
\bibitem {BPV} W.Barth, C.Peters, A.Van de Ven, {\em Compact complex surfaces.}
        Ergebnisse der Mathematik und ihrer Grenzgebiete (3). Springer-Verlag, Berlin,(1984).
\bibitem {B} G.Borrelli, {\em On regular surfaces of general type with $p_g=2$ and
        non birational bicanonical map. }in Algebraic Geometry, A volume in memory of Paolo Francia,
        Beltrametti and alt. (eds.), De Gruyter (2002), 65-78.
\bibitem{C} F.Catanese, {\em Surfaces with $p_g=q=1$ and their period mapping,} in Algebraic Geometry, Lect.
        Notes in Math., $\mathbf{732}$ (1979), 1-26.
\bibitem {CD} F.Catanese,O.Debarre, {\em Surfaces with $K^2=2, p_g=1,q=0$. }J.Reine Angew. Math., 
        $\mathbf{395}$ (1989), 1-55.
\bibitem {CFM} C.Ciliberto, P.Francia, M.Mendes Lopes, {\em Remarks on the bicanonical map for surfaces of
        general type. }Math.Z. $\mathbf {224}$ (1997), 137-166.
\bibitem{CM} C.Ciliberto, M.Mendes Lopes, {\em On regular surfaces of general type with $p_g(S)=3$ and
        non birational bicanonical map. }J.Math. Kyoto Univ., $\mathbf{40}$ (2000), 79-117.
\bibitem {DV} P.Du Val, {\em On surfaces whose canonical system is hyperelliptic}, Canadian J. of Math., 
        $\mathbf 4$ (1952), 204-221. 
\bibitem {F} P.Francia, {\em On the base points of the bicanonical system}, Symposia Mathematica I.N.D.AM.,
        vol. XXXII, Academic Press (1991), 141-150.
\bibitem{H} E.Horikawa, {\em On deformation of quintic surfaces.} Inv.Math., $\mathbf{31}$(1975), 43-85. 
\bibitem {M} D.Morrison, {\em On the moduli of Todorov surfaces.} in Algebraic Geometry and Commutative algebra in 
          honor of Masayoshi NAGATA, (1987), 313-355.
\bibitem{MP1} M.Mendes Lopes, R.Pardini, {\em Enriques surfaces with eight nodes.} Math.Z., $\mathbf{241}$,
             (2002), 673-683.
\bibitem{MP2} M.Mendes Lopes, R.Pardini, {\em The bicanonical map of surfaces with $p\sb g=0$ and $K\sp 2\geq 7$. II.}
              Bull.LononMath.Soc., $\mathbf{35}$, (2003), 337-343.
\bibitem{MP3} M.Mendes Lopes, R.Pardini, {\em The classification of surfaces with $p\sb g=0,K\sp 2=6$ and
              non birational bicanonical map.} (math.AG/0301138).
\bibitem{P} R.Pardini, {\em The classification of double planes of general type with $K\sp 2=8$ and $p\sb g=0$.}
            J.Algebra, $\mathbf{259}$, (2003), 95-118.
\bibitem{Po} F.Polizzi, {\em Surfaces of general type with $p_g=q=1, K^2=8$ and bicanonical map of degree 2.} 
             (math. AG/0311508).
\bibitem {R} I.Reider, {\em Vector bundles of rank 2 and linear systems on algebraic surfaces. }Ann. of Math.,
        $\mathbf{127}$ (1988). 
\bibitem {X1} G.Xiao, {\em Finitude de l'application canonique des surfaces de type general.}
        Bull.Soc.Math.France, $\mathbf{113}$ (1985), 32-51.
\bibitem {X2} G.Xiao, {\em Degree of the bicanonical map of a surface of general type}, Amer. J. of Math.,
        $\mathbf {112}$ (5) (1990), 713-737.

}
\end{thebibliography}
\end{document}